\documentclass{article}
\usepackage{tcolorbox}
\usepackage{amsfonts,amscd,amssymb,amsmath,amsthm,enumerate}
\numberwithin{equation}{section}
\usepackage{array,mathtools,mathrsfs,bm}
\usepackage{threeparttable}
\usepackage{booktabs,caption,longtable,multicol,multirow}
\usepackage{subcaption,lscape}
\usepackage{makecell}
\usepackage{graphicx}
\usepackage[colorlinks,citecolor=purple,linkcolor=blue]{hyperref}
\usepackage[margin=1.5in]{geometry}
\usepackage[numbers]{natbib}
\usepackage{authblk}
\usepackage{algorithm,algorithmic}
\usepackage{algorithm}
\usepackage{pifont}
\usepackage{appendix}
\usepackage{multirow}
\usepackage{tabularx}
\usepackage{makecell}
\usepackage{enumitem}

\usepackage{eqparbox}

\usepackage{tikz}
\definecolor{ao(english)}{rgb}{0.0, 0.5, 0.0}
\definecolor{cadmiumgreen}{rgb}{0.0, 0.42, 0.24}
\definecolor{darkpastelgreen}{rgb}{0.01, 0.75, 0.24}

\usepackage{adjustbox}
\tikzset{%
  materia/.style={draw, fill=blue!20, text width=6.0em, text centered, minimum height=1.5em,drop shadow},
  etape/.style={materia, text width=8em, minimum width=10em, minimum height=3em, rounded corners, drop shadow},
  texto/.style={above, text width=6em, text centered},
  linepart/.style={draw, thick, color=black!50, -LaTeX, dashed},
  line/.style={draw, thick, color=black!50, -LaTeX},
  ur/.style={draw, text centered, minimum height=0.01em},
  back group/.style={fill=white!20,rounded corners, draw=black!50, dashed, inner xsep=15pt, inner ysep=10pt},
}

\tikzstyle{matheq} = [node distance=8.75cm, text width=21em, minimum width=1cm,
minimum height=2em, text centered,color=black!50]

\newtheorem{theorem}{Theorem}[section]

\newtheorem{lemma}{Lemma}[section]
\newtheorem{proposition}{Proposition}[section]

\newtheorem{definition}{Definition}[section]

\newtheorem{assumption}{Assumption}[section]

\newcommand{\x}{\mathbf{x}}
\newcommand{\y}{\mathbf{y}}
\newcommand{\z}{\mathbf{z}}
\newcommand{\p}{\mathbf{p}}

\newcommand{\w}{\mathbf{w}}
\newcommand{\I}{\mathbf{I}}

\title{A Discretization Approach for Bilevel Optimization with Low-Dimensional and Non-Convex Lower-Level}
\author{Xiaotian Jiang\thanks{Department of Industrial and Systems Engineering, University of Minnesota, MN. (jian0851@umn.edu). Equal contribution.}
\quad Ioannis Tsaknakis\thanks{Department of Electrical and Computer Engineering, University of Minnesota, MN. (tsakn001@umn.edu). Equal contribution.}
\quad Prashant Khanduri\thanks{Department of Computer Science, Wayne State University, MI. (khanduri.prashant@wayne.edu).}
\quad Mingyi Hong\thanks{Department of Electrical and Computer Engineering, University of Minnesota, MN. ( mhong@umn.edu).}
}

\begin{document}
\maketitle
\begin{abstract}
Bilevel optimization (BLO) problem, where two optimization problems (referred to as upper- and lower-level problems) are coupled hierarchically,  has wide applications in areas such as machine learning and operations research.  Recently, {many first-order algorithms} have been developed for solving bilevel problems with strongly convex and/or unconstrained lower-level problems{; this special structure of the lower-level problem is needed to ensure the tractability} of gradient computation (among other reasons). In this work, we deal with a class of more challenging BLO problems where the lower-level problem is non-convex and constrained. We propose a novel approach that approximates {the value function of the lower-level} problem by first sampling a set of feasible solutions and then constructing an equivalent convex optimization problem. {This convexified value function} is then used to construct a penalty function for the original BLO problem. 
We analyze the properties of the original BLO problem and the newly constructed penalized problem by characterizing the relation between their KKT points, as well as the local and global minima {of the two problems. We then} develop a gradient descent-based algorithm to solve the reformulated problem, and establish its finite-time convergence guarantees. Finally, we conduct numerical experiments to {corroborate the theoretical performance of the proposed algorithm.} 
\end{abstract}

\section{Introduction}
Bilevel optimization is an {important} tool in modeling interaction between two parties, often with {different} objectives.
This is the case in a wide range of applications, including robust learning \cite{robey2024adversarial,pmlr-v162-zhang22ak}, signal processing \cite{Sun_2022,park2022predicting,cbd1b69037ed4c60831a4683cd0597e5}, meta-learning \cite{DBLP:journals/corr/abs-1909-04630,snell2017prototypical}, hyperparameter optimization \cite{pedregosa2022hyperparameter,franceschi2018bilevel}, and reinforcement learning \cite{hong2020two}.
In this work, we aim to solve such problems; more precisely, we focus on the following \textit{optimistic bilevel optimization} problem,
\begin{align}\label{eq:BLO}\tag{BLO}
&\min\limits_{\x \in \mathcal{X},\z \in \mathcal{R}^m}   f(\x,\z)  \\
& \text{s.t. } \z \in \arg\min\limits_{\y \in \mathcal{Y}} g(\x,\y),\nonumber 
\end{align}
where $f:\mathcal{X}\times\mathcal{Y} \to \mathbb{R}$ and $g:\mathcal{X}\times\mathcal{Y} \to \mathbb{R}$ are continuously differentiable, potentially non-convex functions, and $\mathcal{X} \subseteq \mathbb{R}^{n}$ and $\mathcal{Y} \subseteq \mathbb{R}^{m}$ are convex and compact sets. The above problem \eqref{eq:BLO} consists of two optimization tasks in a hierarchy, the upper-level (UL) problem {(i.e., the minimization of $f$ over $\x$ and $\z$)} and the lower-level (LL) problem {(i.e., the minimization of $g$ over $\y$)}; we also refer to the variable $\y$ as the LL variable and to $\x$ as the UL variable, respectively.

The above bilevel problem is challenging, as the LL problem is both \textit{constrained and non-convex}. One class of algorithms that is able to handle the presence of constraints, albeit only in certain special cases, are implicit gradient methods \cite{ghadimi2018approximation,khanduri2021near, hong2020two}. The key idea behind these methods is the computation of the gradient of the implicit function $F(\x)=f(\x,\y^{\ast}(\x))$ where  $\y^{\ast}(\x)=\arg\min_{\y \in \mathcal{Y}}g(\x,\y)$, which can be subsequently used in gradient-descent type iterative schemes. However, the applicability of these methods is restricted to problems where the LL has a unique solution. This is not the case in {the considered \eqref{eq:BLO}} problem as the LL objective is non-convex.

A few approaches are capable of dealing with both the constraints and the non-convexity in the LL problem (at least to a certain extent). The most common of {these approaches are the} value function-based methods \cite{liu2021valuefunctionbased,ye2022bome,gao2024moreau, liu2023value, shen2023penalty,jiang2024primaldualassistedpenaltyapproachbilevel} 
which utilize the notion of a {\it value function} defined as
\begin{align}\label{eq:vf}
    V(\x):=\min\limits_{\y \in \mathcal{Y}} g(\x,\y) \tag{VF}.
\end{align}
{This} value function is then used to replace the LL problem with the following constraint: $g(\x,\y) \leq V(\x)$, effectively transforming \eqref{eq:BLO} to a {single-level} constrained minimization problem. There has been a growing literature on such value function{-based reformulations} of the \eqref{eq:BLO} problem (to be reviewed in Section \ref{sub:related}); however, this approach has a number of critical issues. For example, when the LL function $g(\x,\y)$ is non-convex in $\y$, the value function $V(\x)$ is non-differentiable and non-convex \cite{10502023}, 
% {\red [Reference?]}, 
therefore, it is difficult to design efficient algorithms for it. In fact, for problem \eqref{eq:BLO}, it is not even possible to evaluate $V(\x)$ as that requires obtaining the global minimum of a non-convex problem.

As discussed above, the \eqref{eq:BLO} problem becomes extremely challenging when the lower-level problem is non-convex. The implicit function based approach does not work since the implicit gradient no longer exists, while the { application of value function-based methods becomes difficult due to the non-convexity and the non-differentiability of the value function.}  Therefore, in this work we focus on a special \eqref{eq:BLO} {problem} where $g(\cdot)$ is non-convex in $\y\in\mathbb{R}^m$, but its dimension $m$ is small. We find that for such a specialized non-convex problem, it is possible to leverage its structure to develop an equivalent reformulation and efficient algorithms. Before diving into the details of the proposed ideas, we provide two representative applications. 

\subsection{Motivating Applications}\label{sec:apps}
\hfill \break
{\bf Toll Setting Problem.}
In this problem a weighted graph $G=(V,E)$ is given where the weights denoted by $\{w_{e}\}_{e \in E}$ describe the costs associated with routing traffic through some edge $e \in E$, e.g., the cost of gas in transportation networks. The {goal of} the \textit{toll setting} problem \cite{Didi-Biha2006, kalashnikov2010comparison} is to search for the toll configuration (i.e., position and pricing of toll stations in the network $G$) that {maximizes} the profit for the toll company. This problem is formulated as a constrained bilevel problem of the form \eqref{eq:BLO}. More precisely, for any given toll configuration, the LL task computes the cheapest way to route the traffic through the network. And, in the UL task, the goal is to optimize over the toll configurations that maximize the total amount paid for tolls given that the cheapest routing is selected (i.e., the LL solution). 
\begin{align}
    &\min\limits_{\substack{\x \in \mathbb{R}^{|E|} \\  \x \geq {\bf 0}, {\bf 1}^\top  \x = c } }  -\sum\limits_{(i,j) \in E}  x_{ij} y_{ij}^\ast(x) \\
    &\text{ s.t.} \;\; \y^\ast(\x) \in 
    \begin{cases} \arg\min\limits_{\y \in \mathbb{R}^{|E|}} {\sum\limits_{(i,j) \in E} g(x_{ij}, y_{ij}; w_{ij})}  \\
     \text{s.t.}  
     \begin{cases} {\bf 0} \leq \y \leq  {\bf b_1},        \sum\limits_{(s,i) \in E} y_{si} =   \sum\limits_{(j,t) \in E} y_{jt}=b_2.  \\
       \sum\limits_{(i,j) \in E} y_{ij} - \!\!\!\! \sum\limits_{(j,k) \in E} y_{jk}=0, \forall j \in V \setminus \{s,t\}. % \nonumber \\
       % \sum\limits_{(s,i) \in E} y_{si} =   \sum\limits_{(j,t) \in E} y_{jt}=\bb_2. 
       \end{cases}
       \end{cases}
    \nonumber
\end{align}%
In the above formulation, the variable $\x=(x_{ij})_{(i,j) \in E}$ represents the toll cost across the network and $\y=(y_{ij})_{(i,j) \in E}$ denotes the quantity of traffic associated with each edge. {The potentially non-convex objective $g(x_{ij}, y_{ij}; w_{ij})$ represents the cost of moving $ y_{ij}$ units of flow through edge $(i,j)$ where the tolls and other associated costs are denoted by $x_{ij}$ and $w_{ij}$, respectively.} Note that the total amount of tolls that can be added to the network is limited to $c$ (UL constraint). Moreover, {the LL constraints ensure that} the traffic capacity of each edge is not violated and the amount of traffic entering a node is equal to the amount exiting.  {Note that} the dimension of the LL problem depends on the number of edges in the network. Therefore, in the case of sparse graphs, the dimensionality of the LL variable {will be} relatively small, {i.e., the dimensionality will be at most of the order of $\mathcal{O}(|V|\log|V|)$, where $|V|$ is the number of nodes \cite[Ch. 7]{alma9975595147801701}} 
. 
{These kinds of sparse graphs arise}, for instance, in transportation networks where typically direct routes exist only between neighboring cities.

\noindent
{\bf Ensemble Learning.}
Ensemble learning is a technique that combines two or more machine learning models (base models) in a way that results {in} improved performance. One way to combine the models is to consider a weighting of their outputs, an instance of the methodology which is referred to as \textit{stacking} \cite{wolpert1992stacked}. The weighting is determined by the solution of an optimization problem where, for instance, we minimize the loss of the weighted model on the training set over the weights, i.e., $\min_{\w \in \Delta} \ell\left(\w^\top H(\x);D^{train}\right)$; $\ell(\cdot)$ is the loss function, $\w$ is the {weight} vector across the $m$ (pretrained) base models, $\Delta=\{\w \in \mathbb{R}^{m}|{\bf 0} \leq \w \leq {\bf 1}, {\bf 1}^\top \w=1\}$ is the {constraint} set {for the weight vector}, $H(\cdot) \in \mathbb{R}^{m}$ is the vector of model outputs, $\x$ is the model input and $D^{train}$ is the training set. Nonetheless, this approach might result {in} overfitting. To address this issue we add the regularizer $\rho_{1} \|\w\|_{2}^{2} + \rho_{2} \|\w\|_{1}$ in the objective, where {$\rho_{1}$ and $\rho_{2}$} are the regularization parameters. To select the optimal values of these parameters we consider a second optimization problem (which relies on the first {one}) where we optimize over {$\rho_{1}$ and $\rho_{2}$} such that the performance of the trained model over the validation set is optimized. The combination of these two problems results {in} the following bilevel optimization problem
\begin{align}\label{eq:bilevel_ensemble}
&\min_{\rho_1, \rho_2 \geq 0, \w \in W} \ell \left(\w^\top H(\x);D^{valid}\right) \\
& \text{s.t. } \w \in \arg\min_{\overline{\w} \in W} \ell \left(\overline{\w}^\top H(\x);D^{train}\right) + \rho_{1} \|\overline{\w}\|_{2}^{2} + \rho_{2} \|\overline{\w}\|_{1}.
\end{align}
The UL task corresponds to the learning of the regularization parameters (hyperparameters), while in the LL task, the weight vector is optimized using a regularized loss function. We note that the dimension of the LL problem is equal to the number of base models, which is typically a small number (e.g., 10s of models).

\subsection{Related Work}\label{sub:related}
There are various solution methods developed for unconstrained bilevel optimization problems. One important approach is implicit differentiation, which uses the Implicit Function Theorem \cite{ghadimi2018approximation}
% {[cite]} 
to compute hypergradients by differentiating through the lower-level problem's optimality conditions. {Utilizing the implicit gradient computation, earlier works developed SGD-type algorithms and established finite-time convergence guarantees for solving bilevel problems  \cite{ghadimi2018approximation, hong2020two}.  
Later improvements were proposed, utilizing techniques such as momentum \cite{khanduri2021near}, variance reduction \cite{yang2021provably}, or conjugate gradient approximations \cite{ji2021bilevel} in order to obtain better convergence guarantees. Nonetheless, a major drawback of the implicit gradient approach is the requirement of a {(unconstrained)} strongly-convex LL problem (in $\y$).} {This requirement is necessary to ensure the differentiability of the implicit function. In addition, implicit gradient methods have been developed for problems with linear equality \cite{xiao2023alternating} and inequality constraints \cite{khanduri2023linearly,tsaknakis2022implicit, kornowski2024first} in the LL. However, in the case of inequality constrained problems, strong convexity of the LL does not suffice to ensure differentiability. In this case, either the problem is modified in a way that {ensures differentiability of the \eqref{eq:BLO} problem} \cite{khanduri2023linearly,tsaknakis2022implicit} or it is {solved} directly using non-smooth optimization methods. 
}

On the other hand, value-function-based methods \cite{liu2021valuefunctionbased,ye2022bome,gao2024moreau, liu2023value, shen2023penalty,jiang2024primaldualassistedpenaltyapproachbilevel} are technically able to handle both the constraints and the non-convexity of the LL problem.  
As mentioned above, the value function approach substitutes the LL problem with an inequality constraint, reducing the bilevel problem into a single-level one. The solution of this single-level problem is typically obtained with the use of penalty-based methods, where the value function constraint is incorporated into the objective using barrier or penalty functions. More broadly, studies on the application of penalty methods for bilevel optimization date back to the 1990s \cite{articlewh, 58565, Ye1997ExactPA, Ishizuka1992DoublePM, aiyoshi1984solution, lin2014solving}, typically 
establishing the equivalence between a relaxation of the original bilevel problem and its penalized form.
However, there are two main issues with such value-function-based approaches. First, the value function is often non-smooth and non-convex, even if the original LL problem is convex in $\y$ \cite{10502023}.
% {\red[cite]}. 
Second, when the LL problem is non-convex, the solutions obtained by penalty-based methods may not be feasible points with respect to the original problem in the lower-level sense, i.e., $g(\x,\y)-V(\x)\ge c$ for some large positive $c$.
Specifically, when penalizing the term \( g(\x,\y) - V(\x) \) as in \cite{shen2023penaltybased}, we obtain a solution \( (\widetilde{\x}, \widetilde{\y}) \) to the penalized problem such that $\left\| \nabla_y \left(g(\widetilde{\x}, \widetilde{\y}) - V(\widetilde{\x}) \right) \right\|$
is small. However, without assuming the Polyak-Łojasiewicz (PL) condition or strong convexity in \cite{shen2023penaltybased}, a small gradient norm does not necessarily imply that
$
g(\widetilde{\x}, \widetilde{\y}) - V(\widetilde{\x})
$
is also small. Therefore, the solution \( (\widetilde{\x}, \widetilde{\y}) \) may be infeasible for the original bilevel problem. In addition to the direct penalty value function approach, there also exist some value-function-based methods \cite{lin2014solving, liu2021valuefunctionbased}. However, these methods either lack theoretical convergence analysis \cite{liu2021valuefunctionbased}, or require technical assumptions to ensure algorithmic convergence, such as \cite[Assumption 3.1]{lin2014solving}. In contrast, our algorithm does not rely on any non-standard assumptions for convergence, nor does it require the PL condition on the lower-level problem as in \cite{shen2023penaltybased}.

\subsection{Contributions}

In this work, we develop a new approach for solving a class of bilevel optimization problems {with non-convex and constrained LL problems.} More precisely, the contributions of this work are the following:
\begin{enumerate}[leftmargin = *]
    \item We propose a value-function-based reformulation for {solving} bilevel problems of the form \eqref{eq:BLO} with non-convex LL constraints. The key novelty is the \textit{discretization of the value function}, which results in a reformulated problem that is easy to optimize (e.g., the reformulation is smooth)
    at a price of an approximation error that can directly be controlled by a regularization parameter and the number of {discretization} points.
    \item {We establish the equivalence between the global solutions of the reformulated discretized problem and the $\epsilon$-relaxation of the original bilevel problem, under certain assumptions. In addition, we consider a penalty-based reformulation of the discretized problem and show that these two problems share the same set of local minima, global minima, and KKT points. These results justify the significance of the discretized formulation and support the validity of solving the penalized problem as a surrogate for the discretized formulation.} 
    \item {To efficiently solve the} reformulated problem, {we develop} the {DIscretized Value-FunctIon-BaseD - Bilevel Optimization (DIVIDE-BLO)} algorithm, and establish its finite-time convergence guarantees to reach certain stationary points of the reformulated problem.
    \item We conduct numerical experiments on synthetic and real-world (ensemble learning) problems that showcase the utility of the proposed approach.
\end{enumerate}

\section{Bilevel Discretization Approach}\label{sec:reform}

To begin with, we introduce assumptions that hold throughout the paper.
\begin{assumption}\label{ass:basics}
The following assumptions hold:
\begin{enumerate}[label=(\alph*), leftmargin = *]
\item The functions $f(\x,\y)$ and $g(\x,\y)$ are continuously differentiable.
\item\label{ass:basics:sets} The sets $\mathcal{X}$ and $\mathcal{Y}$ are convex and compact with $\|\y\|\leq D$ for any $\y$ in $\mathcal{Y}$ {and $D>0$ is a constant}.
\item\label{ass:basics:Lip_grad_f} The function $f$ has Lipschitz continuous gradients, i.e., 
$$\|\nabla f(\x_{1},\y_{1}) - \nabla f(\x_{2},\y_{2})\| \leq \overline{L}_{f} \|[\x_{1}; \y_{1}]-[\x_{2}; \y_{2}]\|, \forall \x_{1},\x_{2} \in \mathcal{X}, \y_{1},\y_{2} \in \mathcal{Y}.$$
\item\label{ass:basics:Lip_grad_g} The function $g$ has Lipschitz continuous gradients, i.e., 
$$\|\nabla g(\x_{1},\y_{1}) - \nabla g(\x_{2},\y_{2})\| \leq \overline{L}_{g} \|[\x_{1}; \y_{1}]-[\x_{2}; \y_{2}]\|, \forall \x_{1},\x_{2} \in \mathcal{X}, \y_{1},\y_{2} \in \mathcal{Y}.$$
\item\label{ass:basics:Lip_bound_f} The function $f$ has bounded gradients, i.e. $\|\nabla f(\x,\y)\|\leq L_f, \forall \x \in \mathcal{X}, \y \in \mathcal{Y}$
\item\label{ass:basics:Lip_bound_g} The function $g$ has bounded gradients, i.e. $\|\nabla g(\x,\y)\|\leq L_g, \forall \x \in \mathcal{X}, \y \in \mathcal{Y}$
\item For every $\x \in \mathcal{X}$, the LL problem is feasible, i.e., it has at least one solution.
\end{enumerate}
\end{assumption}
These assumptions are standard in the literature of bilevel optimization, e.g., see \cite{ghadimi2018approximation,hong2020two,chen2021singletimescale}.

\subsection{Approximate Value Function and Its Properties}

The proposed approach is based on the \textit{discretization of the LL problem and its value function}. To be more precise, {we} select $k$ points from the LL {constraint} set $\mathcal{Y}$ {denoted as} $\widetilde{\mathcal{Y}}=\left\{\y^{(1)},\y^{(2)},\ldots,\y^{(k)}\right\} \subseteq \mathcal{Y}$, such that these points provide a covering of radius $r$ of the set {complete set} $\mathcal{Y}$. This implies that for every $\y \in \mathcal{Y}$, there exists an index $j \in [k]$ such that $\|\y-\y^{(j)}\| \leq r$. It can be shown that for a given bounded set $\mathcal{Y} \subseteq \mathbb{R}^{m}$ with $\|\y\| \leq D, \forall \y \in \mathcal{Y}$, {the minimum number of points $k$ {required} to cover the set $\mathcal{Y}$ is at most $\left(2D\sqrt{m}/r\right)^{m}$ with a covering radius $r$ \cite{shalev2014understanding}.} 
{
We assume throughout the paper that we have such a covering of the set $\mathcal{Y}$. To emphasize this point we introduce the following assumption.
\begin{assumption}\label{ass:cover}
    The set $\widetilde{\mathcal{Y}}=\left\{\y^{(1)},\y^{(2)},\ldots,\y^{(k)}\right\} \subseteq \mathcal{Y}$ defines a covering of the set $\mathcal{Y}$ with radius $r$.
\end{assumption}
}
{In practice, we can build the discretized set $\widetilde{\mathcal{Y}}$, either by faithfully dividing the entire space deterministically into a grid (for relatively simple constraints) or by randomly sampling a sufficiently large number of points from the set $\mathcal{Y}$ (for complex constraints). Although for complex constraints sets we cannot provably obtain a covering of a given radius in that way, we notice that, in practice, this approach works satisfactorily, as demonstrated in our experiments (see Sec. \ref{sec:exp_ensemble}).} 

Using the discretization of $\mathcal{Y}$ {with} $\widetilde{\mathcal{Y}}$, we formulate the approximate value function $$\widehat{V}(\x)=\min\limits_{i \in [k]} \left\{ g\left(\x,\y^{(i)}\right)\right\},$$ which involves {solving} a discrete optimization problem. {Note} that the underlying optimization problem of the value function can be {equivalently} written as a continuous convex optimization problem over the simplex constraints {as} $$\min\limits_{{\bf 0} \leq \p \leq {\bf 1}, {\bf 1}^\top \p=1} \left\{ \sum_{i=1}^{k} p_{i}g\left(\x,\y^{(i)}\right)\right\}.$$ Finally, we add a quadratic regularizer $(\lambda/2)\|\p\|^{2}$, for some $\lambda>0$, to ensure that the above minimization problem is strongly convex and therefore has a unique solution, 
\begin{align}\label{eq:DVF}
    \widetilde{V}_{\lambda}(\x)=\min\limits_{\p \in \Delta} \left\{ \sum\limits_{i=1}^{k} p_{i}g\left(\x,\y^{(i)}\right) + \frac{\lambda}{2}\|\p\|^{2}\right\}. \tag{DISC-VF}
\end{align}
where we denote {the unit simplex by} $\Delta=\{\p \in \mathbb{R}^{k}|{\bf 0} \leq \p \leq {\bf 1}, {\bf 1}^\top \p=1\}$.

We would like to stress here that the approximate value function $\widetilde{V}_{\lambda}(\x)$ has some desirable properties as compared to the original one $V(\x)$. {Specifically}, the minimization problem in {the definition of} $V(\x)$ is non-convex, while the one in the discretized one $\widetilde{V}_{\lambda}(\x)$ is an easier (strongly) convex problem over the variable $\p$. The strong convexity of the underlying problem ensures the uniqueness of the solution, which renders the approximate value function $\widetilde{V}_{\lambda}(\x)$ differentiable. {Note} that {when lower-level problem is non-convex in $\y$}, {the original value function} $V(\x)${$=\arg\min_{\y\in\mathcal{Y}}g(\x,\y)$} is {non-differentiable in general} \cite{10502023}.
Naturally, one might introduce the quadratic regularizer $(\bar{\lambda}/2)\|\y\|^{2}, \bar{\lambda}>0$ in the minimization problem of $V(\x)$ {to} make the objective strongly convex, {similar} to our proposed approach. However, in the case of $V(\x)$, $\bar{\lambda}$ has to be sufficiently large, as the objective $g(\x,\y)$ is non-convex. This results in an approximation error that { cannot be} directly controlled. {In contrast, for} the approximate value function $\widetilde{V}_{\lambda}(\x)$ the regularization parameter $\lambda$ can be chosen \textit{arbitrarily small} to maintain the strong convexity of the objective, resulting in improved control over the approximation error. Finally, it is worth mentioning that the discretization of the value function introduces an error, as it is an approximation of the original one $V(\x)$. However, we show that this error becomes arbitrarily small as the number of points $k$ becomes larger and the regularization parameter $\lambda$ becomes smaller.
We formally establish the aforementioned results in the following proposition.
\begin{proposition}\label{pro:dvf_properties}
{Under Assumptions \ref{ass:basics}, \ref{ass:cover}, the following statements} hold:
\begin{enumerate}[label=(\alph*),leftmargin = *]
    \item The approximate value function $\widetilde{V}_{\lambda}(\x)$ is differentiable, and its gradient is given by the following {expression}
    \begin{align}
        \nabla \widetilde{V}_{\lambda}(\x) = \sum\limits_{i=1}^{k} p_{i}^{\ast}\nabla_{\x}g(\x,\y^{(i)})  
    \end{align}
    where $p_{i}^{\ast}$ are the elements of {vector} $\p^{\ast} = \arg\min\limits_{\p \in \Delta} \left\{ \sum\limits_{i=1}^{k} p_{i}g(\x,\y^{(i)}) + \frac{\lambda}{2}\|\p\|^{2}\right\}$.
    \item\label{pro:dvf_grad:error} The approximate value function $\widetilde{V}_{\lambda}(\x)$ approximates $V(\x)$ with the following approximation error $$\big\lvert V(\x)-\widetilde{V}_{\lambda}(\x)\big\lvert\leq \frac{2L_gD\sqrt{m}}{k^{1/m}}+\frac{\lambda}{2}, \; \forall \x \in \mathcal{X}.$$
    \item {The approximate value function} $\widetilde{V}_{\lambda}(\x)$ has Lipschitz continuous {gradient, i.e.,} 
    $$\|\nabla \widetilde{V}_{\lambda}(\x_{1}) - \nabla \widetilde{V}_{\lambda}(\x_{2})\| \leq L(\lambda,k) \|\x_{1}-\x_{2}\|, \forall \x_{1},\x_{2} \in \mathcal{X},$$
    where {$L(\lambda,k)=\frac{L_g^2k}{\lambda}+\overline{L}_g$} {denotes} the Lipschitz constant. 
\end{enumerate}
\end{proposition}
\begin{proof}
    {See} Appendix \ref{app:proofs:reform}.
\end{proof}

\subsection{Discretized Bilevel Problem Reformulation}

To begin with our approach, let us first provide a (relaxed) value function reformulation of \eqref{eq:BLO} 
\begin{align}\label{eq:blo_vf}\tag{BLO-VF}
&\min\limits_{\x \in \mathcal{X},\y \in \mathcal{Y}}   f(\x,\y)  \\
& \text{s.t. } g(\x,\y) \leq V(\x)  + \epsilon. \nonumber
\end{align}
where $\epsilon>0$ is some small constant introduced to ensure {strict feasibility of the LL constraint, i.e.,} $\epsilon$ ensures that there are points on which the value function inequality constraint is strictly feasible. {This is a standard practice in the relevant literature, e.g., see \cite{lin2014solving, shen2023penaltybased}}. 
We perform a \textit{discretization of the LL problem} of \eqref{eq:blo_vf} over the set $\widetilde{\mathcal{Y}}$. This involves using value function $\widetilde{V}_{\lambda}(\x)$ in place of $V(\x)$, as described above, and the reparameterization of the reformulation \eqref{eq:blo_vf}, where the variable $\y$ is {replaced} by a variable $\p$ that belongs to the {$k$-dimensional} unit simplex. The new variable $\p$ allows us to express $g(\x,\y)$ using only its values on the discrete set of points $\widetilde{\mathcal{Y}}$ as $\sum_{i=1}^{k}p_ig(\x,\y^{(i)})$. 
The {discretized formulation} is provided below:
\begin{align}\label{eq:disc_reform}\tag{BLO-DISC}  % \label{eq:disc_reform_1}
    &\min\limits_{{\x \in \mathcal{X}},\p} f\left(\x,\sum\limits_{i=1}^{k} p_{i}\y^{(i)}\right)  \\
    & \text{s.t. }  \sum\limits_{i=1}^{k} p_{i}g\left(\x,\y^{(i)}\right) - \widetilde{V}_{\lambda}(\x) \leq \epsilon \nonumber\\
    & \quad\quad {\bf 0} \leq \p \leq {\bf 1}, {\bf 1}^\top \p=1 \nonumber
\end{align}
{where the approximate value function $\widetilde{V}_{\lambda}(\x)$ is defined in equation \eqref{eq:DVF}.} 

To simplify notation, we define $\widetilde{f}(\x,\p):=f(\x,\sum_{i=1}^{k} p_{i}\y^{(i)})$. The discretized problem \eqref{eq:disc_reform} {stated above} is a difficult constrained minimization task. The major source of the difficulty comes from the presence of the non-convex value function $\widetilde{V}_{\lambda}(\x)$ in the constraint. 
To proceed,  we adopt a penalty-based approach to penalize the value function constraint, and consider the following problem: 
\begin{align}\label{eq:pen_reform}\tag{BLO-PEN}
    &\min\limits_{{ \x \in \mathcal{X}},\p} F_{\gamma}(\x,\p):= f\left(\x,\sum\limits_{i=1}^{k} p_{i}\y^{(i)}\right) + \gamma \left(\sum\limits_{i=1}^{k} p_{i}g\left(\x,\y^{(i)}\right) - \widetilde{V}_{\lambda}(\x) \right) \\
    & \text{s.t. } {\bf 0} \leq \p \leq {\bf 1}, {\bf 1}^\top \p=1 \nonumber. 
\end{align}

\section{Equivalence Between Different Problems}\label{sec:EquivalenceBetweenDifferentProblems}
So far, we have introduced four different problems: the original problem \eqref{eq:BLO}, the value-function-based relaxation \eqref{eq:blo_vf}, the discretized reformulation \eqref{eq:disc_reform}, and the penalty-based discretized reformulation \eqref{eq:pen_reform}. Next we will study the relationship  among \eqref{eq:blo_vf}, \eqref{eq:disc_reform}, and \eqref{eq:pen_reform}. 
The reason for excluding \eqref{eq:BLO} from our subsequent discussion is twofold. First, without additional assumptions, it is extremely difficult to establish any meaningful connection between the original problem and the relaxed formulation \eqref{eq:blo_vf}. In particular, when the lower-level problem is nonconvex, even a slight relaxation may dramatically enlarge the feasible set, making it intractable to analyze the equivalence between the two. Second, the relaxed problem itself is often taken as the primary object of study in the literature—for instance, see \cite{shen2023penaltybased,lin2014solving}.

\begin{figure}[htbp]
\centering

\scalebox{0.785}{

\tikzset{every picture/.style={line width=0.75pt}} %set default line width to 0.75pt        

\begin{tikzpicture}[x=0.75pt,y=0.75pt,yscale=-1,xscale=1]

%Straight Lines [id:da02741071658274752] 
\draw  [dash pattern={on 4.5pt off 4.5pt}]  (50,170) -- (77.2,170) ;
\draw [shift={(79.2,170)}, rotate = 180] [color={rgb, 255:red, 0; green, 0; blue, 0 }  ][line width=0.75]    (10.93,-3.29) .. controls (6.95,-1.4) and (3.31,-0.3) .. (0,0) .. controls (3.31,0.3) and (6.95,1.4) .. (10.93,3.29)   ;
%Straight Lines [id:da9737125734865333] 
\draw    (179,169.5) -- (262.2,169.5)(179,172.5) -- (262.2,172.5) ;
\draw [shift={(270.2,171)}, rotate = 180] [color={rgb, 255:red, 0; green, 0; blue, 0 }  ][line width=0.75]    (10.93,-3.29) .. controls (6.95,-1.4) and (3.31,-0.3) .. (0,0) .. controls (3.31,0.3) and (6.95,1.4) .. (10.93,3.29)   ;
\draw [shift={(171,171)}, rotate = 0] [color={rgb, 255:red, 0; green, 0; blue, 0 }  ][line width=0.75]    (10.93,-3.29) .. controls (6.95,-1.4) and (3.31,-0.3) .. (0,0) .. controls (3.31,0.3) and (6.95,1.4) .. (10.93,3.29)   ;
%Straight Lines [id:da40790058593317946] 
\draw    (389,169.5) -- (511.2,169.5)(389,172.5) -- (511.2,172.5) ;
\draw [shift={(519.2,171)}, rotate = 180] [color={rgb, 255:red, 0; green, 0; blue, 0 }  ][line width=0.75]    (10.93,-3.29) .. controls (6.95,-1.4) and (3.31,-0.3) .. (0,0) .. controls (3.31,0.3) and (6.95,1.4) .. (10.93,3.29)   ;
\draw [shift={(381,171)}, rotate = 0] [color={rgb, 255:red, 0; green, 0; blue, 0 }  ][line width=0.75]    (10.93,-3.29) .. controls (6.95,-1.4) and (3.31,-0.3) .. (0,0) .. controls (3.31,0.3) and (6.95,1.4) .. (10.93,3.29)   ;

% Text Node
\draw (2,162.4) node [anchor=north west][inner sep=0.75pt]    {$(BLO)$};
% Text Node
\draw (51,153) node [anchor=north west][inner sep=0.75pt]  [font=\footnotesize] [align=left] {relax};
% Text Node
\draw (82,162.4) node [anchor=north west][inner sep=0.75pt]    {$(BLO-VF)$};
% Text Node
\draw (172,153) node [anchor=north west][inner sep=0.75pt]  [font=\footnotesize] [align=left] {$\displaystyle \delta -$global solution};
% Text Node
\draw (172,183) node [anchor=north west][inner sep=0.75pt]  [font=\footnotesize] [align=center] {under Assumption \ref{ass:reform}};
% Text Node
\draw (272,162.4) node [anchor=north west][inner sep=0.75pt]    {$(BLO-DISC)$};
% Text Node
\draw (522,162.4) node [anchor=north west][inner sep=0.75pt]    {$(BLO-PEN)$};
% Text Node
\draw (382,153) node [anchor=north west][inner sep=0.75pt]  [font=\footnotesize] [align=left] {global/local/KKT solution};

\end{tikzpicture}
}
\caption{Equivalence between different problems.}
\label{fig:blo-equiv}
\end{figure}
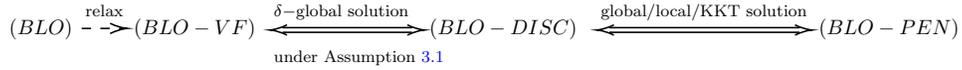

More specifically, we will establish the $\delta$-global 
solution equivalence between \eqref{eq:blo_vf} and \eqref{eq:disc_reform} in Section \ref{sec:equiv2and3}, and the global/ local/ KKT solution equivalence between \eqref{eq:blo_vf} and \eqref{eq:disc_reform} in Section \ref{sec:equiv3and4}; see Figure \ref{fig:blo-equiv} for an illustration of the relationships between different problems.

\subsection{Equivalence Between \eqref{eq:blo_vf} and \eqref{eq:disc_reform}}\label{sec:equiv2and3}

In this section, we investigate the relation between the discretized bilevel problem reformulation (\ref{eq:disc_reform}) and the relaxed problem (\ref{eq:blo_vf}). Notably, without imposing any additional assumptions, the relaxed problem (\ref{eq:blo_vf}) is even more challenging than solving single-level non-convex optimization problems. In single-level non-convex optimization, we just need to find one solution whose function value is close to the global minimum. In contrast, the relaxed problem (\ref{eq:blo_vf}) requires us to study the set of all points whose values are close to the minimum. Specifically, for any $\x$, to evaluate the following hyperfunction value
$$\min_{\y\in\mathcal{Y}}f(\x,\y)\quad s.t.\quad g(\x,\y)-V(\x)\leq\epsilon,$$
we need to find the point that minimizes $f(\x,\y)$ over the level set $\{\y\in\mathcal{Y}:g(\x,\y)-V(\x)\leq\epsilon\}$. Unfortunately, even when employing a sampling method that exhaustively covers a sufficiently large number of points in the lower-level domain, there is no guarantee of finding this kind of point. 

\begin{figure}[htbp]
\centering

\begin{minipage}{.45\textwidth}
\centering
\includegraphics[width=\linewidth]{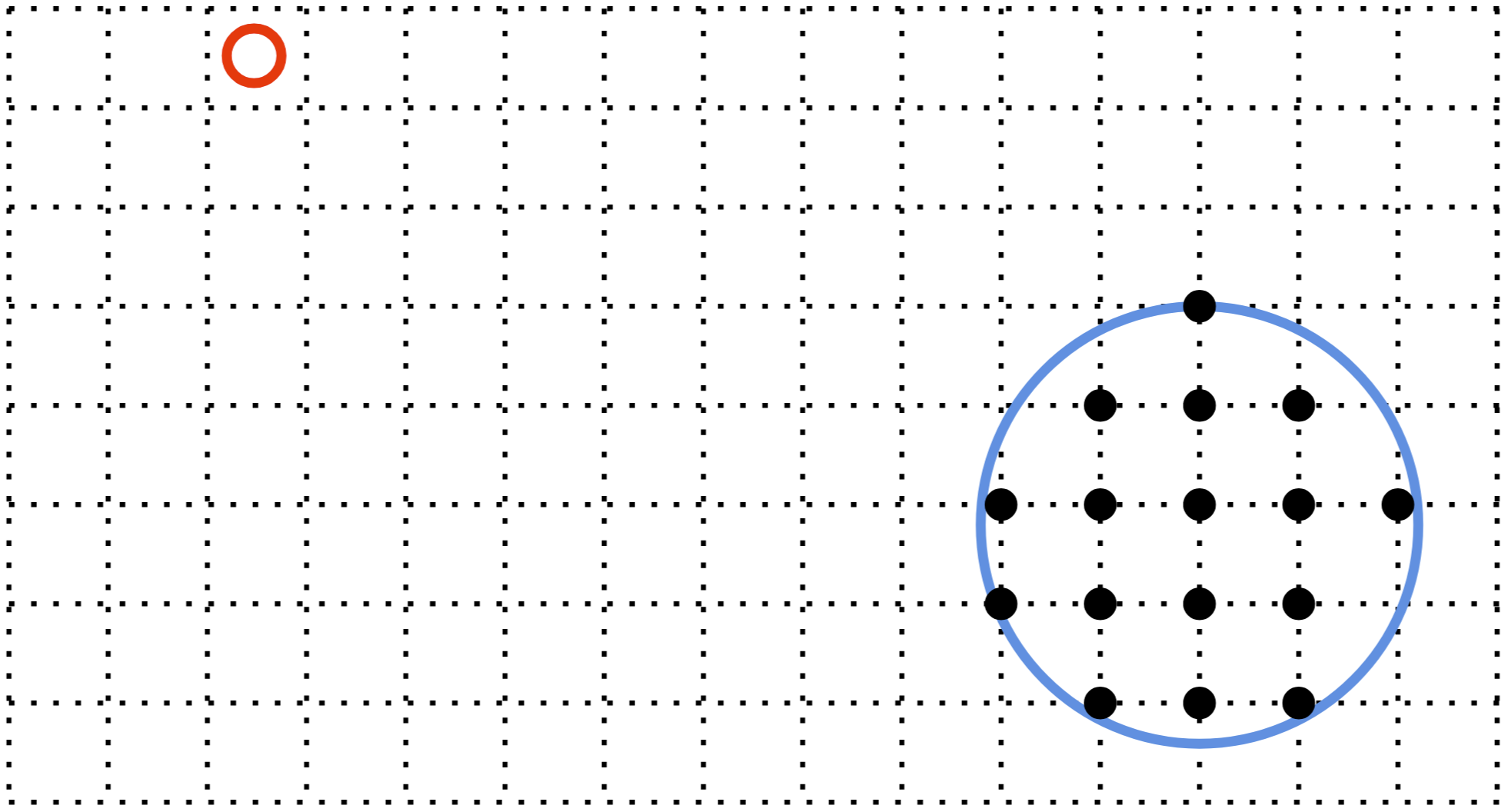}
\caption{This figure depicts the constraint set $\mathcal{Y}$ and its discretization; the sample points are the grid intersections. The red region in the top-left corner and the blue region in the bottom-right corner are the feasible regions within $\mathcal{Y}$, i.e., subsets of $\left\{ \y \in \mathcal{Y}: g(\x, \y) - V(\x) \leq \epsilon \right\}.$ The feasible region in the bottom-right corner is captured by the sample points (i.e., the highlighted vertices of the grid), but the feasible region in the top-left corner is too small and is not captured by any sample point.}
\label{fig:feasibleset}
\end{minipage}%
\hspace{5mm}
\begin{minipage}{.45\textwidth}
\centering
\includegraphics[width=\linewidth]{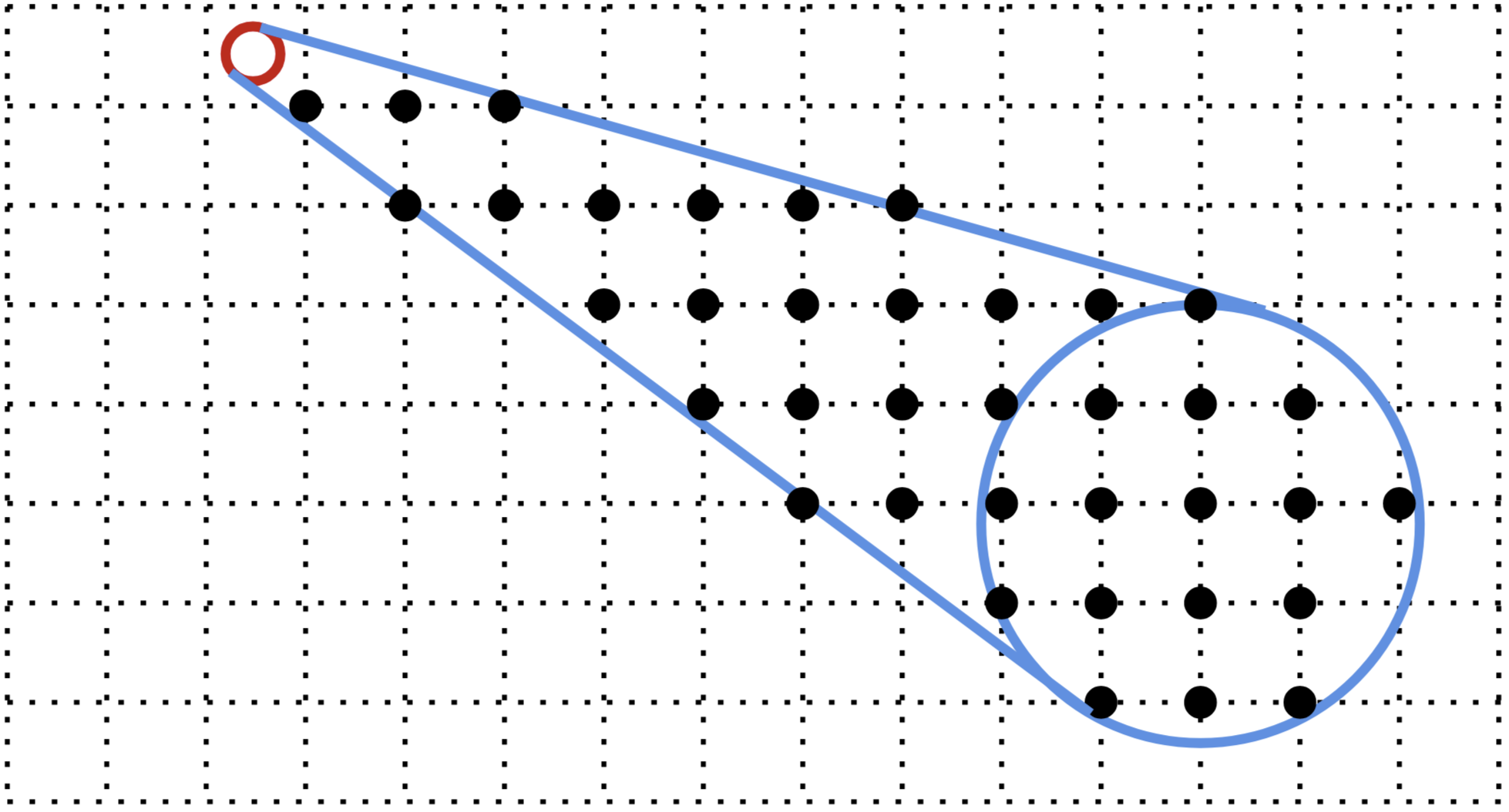}
\caption{This figure depicts the constraint set $\mathcal{Y}$ and its discretization; the sample points are the grid intersections. When $\epsilon$ is sufficiently small, the level set $\left\{ \y \in \mathcal{Y} : g(\x, \y) - V(\x) \leq \epsilon \right\}$ is convex. Therefore, the convex hull of the top-left and bottom-right regions (delineated by the blue lines) is part of the level set, and our sample points are able to capture the entire level set.\\
\\
\\}
\label{fig:convexhull}
\end{minipage}

\end{figure}

For example, no matter how finely we select our sample points, we might encounter the following situation: there exists one point $\y$ such that $g(\x,\y)-V(\x)$ is just slightly less than $\epsilon$, but a slight deviation from $\y$ immediately makes $g(\x,\y)-V(\x)$ exceed $\epsilon$. {That is, the feasible region in the neighborhood of $\y$ is very small.} 
As a result, the sample points might fail to capture such a $\y$, which may yield a very small value of $f(\x,\y)$. This situation is depicted in Figure \ref{fig:feasibleset}, where the grid intersections are our sample points. These sample points successfully capture the circular feasible region in the bottom right corner. However, the region in the top left is very small, so the sample points happen to miss it. It is possible that the value of $\y$ within that missed region actually minimizes $f(\x,\y)$.

To make the analysis of the equivalence between \eqref{eq:blo_vf} and \eqref{eq:disc_reform} tractable, we need the following assumptions.

\begin{assumption}\label{ass:reform}
There exists positive constants $\epsilon_0$, $C$ and $\beta$ such that the following conditions hold:
\begin{enumerate}[leftmargin = *]
    \item For any $\x \in \mathcal{X}$ and $\epsilon\leq\epsilon_0$, the level set $\{\y\in\mathcal{Y}:g(\x,\y)-V(\x)\leq\epsilon\}$ is convex, where $V(\x)=\min\limits_{\y \in \mathcal{Y}} g(\x,\y)$ is the value function;\label{ass:reform(1)}
    \item For any $\y\in\{\y\in\mathcal{Y}:g(\x,\y)-V(\x)\leq\epsilon_0\}$, we have $\text{dist}(\y,\arg\min\limits_{\y\in\mathcal{Y}}g(\x,\y))\leq C\cdot |g(\x,\y)-V(\x)|^\beta$.\label{ass:reform(2)}
\end{enumerate} 
\end{assumption}

It is important to note that Assumption \ref{ass:reform} is only required in this section, specifically for analyzing the equivalence between \eqref{eq:blo_vf} and \eqref{eq:disc_reform}. These assumptions are not needed for analyzing the equivalence between \eqref{eq:disc_reform} and \eqref{eq:pen_reform}, nor for solving \eqref{eq:pen_reform}.

{The first condition in} Assumption \ref{ass:reform} ensures that the counterexample in Figure \ref{fig:feasibleset} does not occur. This is because the level set is convex in this case and our sample points can capture almost the entire feasible region (see Figure \ref{fig:convexhull}). Moreover, the assumption implies that if several sampled points nearly achieve the lowest {upper-level} function value, then any point formed by a convex combination of these points will also have a {small} function value. As a result, it guarantees that the $(\x,\y)$ corresponding to a feasible solution $(\x,\p)$ of problem (\ref{eq:disc_reform}) is also feasible for (\ref{eq:blo_vf}), as demonstrated in the following Lemma \ref{lem:ll_bound}. 
{The second condition in} Assumption \ref{ass:reform} allows us to more precisely estimate the distance between $\y$ and the optimal set $\y^\ast(\x)$ based on the difference between the function value $g(\x,\y)$ and $V(\x)$. As a result, it enables us to analyze the $\delta$-global solution equivalence between \eqref{eq:blo_vf} and \eqref{eq:disc_reform}, where a $\delta$-global solution is defined in Definition \ref{def:deltasolution}.

\begin{lemma}\label{lem:ll_bound}
    Suppose Assumptions \ref{ass:basics}, \ref{ass:cover}, \ref{ass:reform} hold, $$k\geq\max\{(2L_gD\sqrt{m}/\epsilon_0)^m,(4L_gD\sqrt{m}/\epsilon_1)^m\},$$ and $\lambda\leq \epsilon_1\leq \min\{1,\epsilon_0^2\}$. If $\sum_{i=1}^kp_ig(\x,\y^{(i)})-\widetilde{V}_{\lambda}(x)\leq\epsilon_1$, then 
    $$g\left(\x,\sum_{i=1}^kp_i\y^{(i)}\right)-V(\x)\leq \left(1+4D\sqrt{m}L_g \right)\sqrt{\epsilon_1}.$$
\end{lemma}
\begin{proof}
    See Appendix \ref{app:proofs:Equivalence}.
\end{proof}
    The main idea of the proof of Lemma \ref{lem:ll_bound} is as follows: Assuming that $\epsilon_1$ is small enough, we divide the sample points into two categories. The \textit{first category} consists of sample points for which the difference between the function value $g$ and the minimum value is very small, i.e. those satisfying $g(\x,\y^{(i)})-V(\x)\leq\epsilon$, where $\epsilon$ is a small constant related to $\epsilon_1$. The remaining sample points form the \textit{second category}. Note that our goal is to estimate the difference between the function value corresponding to the weighted combination of these sample points and $V(\x)$, i.e. $g(\x,\sum_{i=1}^kp_i\y^{(i)})-V(\x)$. If only the sample points in the first category have positive weights $p_i$, then by Assumption \ref{ass:reform}.\ref{ass:reform(1)} (i.e. the convexity of the level set), we know that the function value of their linear combination will also be close to $V(\x)$. If some sample points in the second category also have positive weights, then as long as the sum of these weights is small, the Lipschitz continuity of $g(\x,\y)$ ensures that their influence on the overall weighted combination is negligible. Therefore, the key difficulty of this argument lies in choosing an appropriate $\epsilon$ and evaluating the influence of the second category of sample points on the overall weighted combination.

\begin{definition}[Approximate global solutions of \eqref{eq:blo_vf} and \eqref{eq:disc_reform}.]\label{def:deltasolution}
\begin{enumerate}[leftmargin = *]
    \item We say that $(\x, \y)$ is a $\delta$-global solution of the relaxed problem~\eqref{eq:blo_vf} with relaxation coefficient $\epsilon$, if it satisfies the relaxed feasibility condition $g(\x, \y) - V(\x) \leq \epsilon$, and for any $(\x', \y')$ such that $g(\x', \y') - V(\x') \leq \epsilon$, it holds that $f(\x, \y) \leq f(\x', \y') + \delta$. 
    \item We say that $(\x, \p)$ is a $\delta$-global solution of the discretized bilevel problem~\eqref{eq:disc_reform} with relaxation coefficient~$\epsilon$, if it satisfies the relaxed feasibility condition $\sum_{i=1}^k p_i g(\x, \y^{(i)}) - \widetilde{V}_\lambda(\x) \leq \epsilon$, and for any $(\x', \p')$ such that $\sum_{i=1}^k p'_i g(\x', \y^{(i)}) - \widetilde{V}_\lambda(\x') \leq \epsilon$, it holds that $\widetilde{f}(\x, \p) \leq \widetilde{f}(\x', \p') + \delta$. 
\end{enumerate}
\end{definition}

We now formally present the equivalence between \eqref{eq:blo_vf} and \eqref{eq:disc_reform} with respect to $\delta$-global solutions.

\begin{theorem}\label{thm:global}
    Suppose Assumptions \ref{ass:basics}, \ref{ass:cover}, \ref{ass:reform} hold, $$k\geq\max\{(2L_gD\sqrt{m}/\epsilon_0)^m,(4L_gD\sqrt{m}/\epsilon_1)^m\},$$ and $0<\lambda\leq\epsilon_1\leq\min\{1,\epsilon_0^2\}$. We have the following equivalence:
    \begin{enumerate}[leftmargin = *]
        \item \label{thm:global:part1}If $(\widetilde{\x},\widetilde{\p})$ is a global solution of \eqref{eq:disc_reform} with relaxation coefficient $\epsilon_1\leq[\epsilon_0/(1+4DL_g)]^2$, then $(\widetilde{\x},\widetilde{\y})=(\widetilde{\x},\sum_{i=1}^k\widetilde{p}_i\y^{(i)})$ is also a $\delta$-global solution of \eqref{eq:blo_vf} with relaxation coefficient $(1+4DL_g)\sqrt{\epsilon_1}$. Here $\delta$ is defined as follows
        \begin{align*}
            \delta:=L_f\cdot\left\{C\left[\left(1+4DL_g\right)\sqrt{\epsilon_1}\right]^\beta+\frac{\epsilon_1}{2L_g}\right\}.
        \end{align*}
        \item If $(\widetilde{\x},\widetilde{\y})$ is a global solution of \eqref{eq:blo_vf} with relaxation coefficient $\epsilon_1\leq[2\epsilon_0/(5+20DL_g)]^2$, then there exists a point $(\widetilde{\x},\widetilde{\p})$ such that $\|(\widetilde{\x},\widetilde{\y})-(\widetilde{\x},\sum_{i=1}^k\widetilde{p}_i\y^{(i)})\|\leq \epsilon_1/(2L_g)$, and $(\widetilde{\x},\widetilde{\p})$ is a $\delta$-global solution of (\ref{eq:disc_reform}) with relaxation coefficient $5\epsilon_1/2$. Here $\delta$ is defined as follows
        \begin{align*}
            \delta:=L_f\cdot\left\{C\left[\frac{5}{2}\left(1+4DL_g\right)\sqrt{\epsilon_1}\right]^\beta+\frac{\epsilon_1}{2L_g}\right\}.
        \end{align*}
    \end{enumerate}
\end{theorem}
\begin{proof}
    See Appendix \ref{app:proofs:Equivalence}.
\end{proof}

For the first part of Theorem \ref{thm:global}, we use a proof by contradiction. If the $(\widetilde{\x},\widetilde{\y})$ induced by $(\widetilde{\x},\widetilde{\p})$ is not a $\delta$-{global} solution, i.e. there exists another $(\x',\y')$ such that $f(\x',\y')<f(\widetilde{\x},\widetilde{\y})-\delta$, then we can find a sample point that is close to $\y'$, from which we can reconstruct a feasible pair $(\x',\p')$ such that $\widetilde{f}(\x',\p')<\widetilde{f}(\widetilde{\x},\widetilde{\p})$, contradicting the assumption that $(\widetilde{\x},\widetilde{\p})$ is a minimizer. For the second part, we note that $(\widetilde{\x},\widetilde{\y})$ may not correspond to any $(\widetilde{\x},\widetilde{\p})$. Therefore, we aim to show that there exists a pair $(\widetilde{\x},\widetilde{\p})$ such that the point it induces is close to $(\widetilde{\x},\widetilde{\y})$, and $(\widetilde{\x},\widetilde{\p})$ is a $\delta$-global solution. We again proceed by contradiction: if not, then there exists a pair $(\x',\p')$ such that $\widetilde{f}(\x',\p')<\widetilde{f}(\widetilde{\x},\widetilde{\p})-\delta$. We can prove that $(\x',\y')$ induced by $(\x',\p')$ is still feasible, and satisfies $f(\x',\y')<f(\widetilde{\x},\widetilde{\y})$, which contradicts the assumption that $(\widetilde{\x},\widetilde{\y})$ is a minimizer.

\subsection{Equivalence Between \eqref{eq:disc_reform} and \eqref{eq:pen_reform}}\label{sec:equiv3and4}
Next, we {establish} the connection between the optimal points of \eqref{eq:disc_reform} and \eqref{eq:pen_reform}. We {first state} a preliminary lemma.
\begin{lemma}\label{lem:KKT1}
Suppose Assumptions \ref{ass:basics},\ref{ass:cover} hold and $\lambda \leq \frac{L_f D}{\gamma}$. If $(\x^{\ast},\p^{\ast})$ is a KKT point of \eqref{eq:pen_reform}, then
    \begin{align}
        \sum\limits_{i=1}^{k} p_{i}^{\ast}g\left(\x^{\ast},\y^{(i)}\right) - \widetilde{V}_{\lambda}(\x^{\ast})\leq\frac{5L_f D}{2\gamma}.
    \end{align}
\end{lemma}
\begin{proof}
    See Appendix \ref{app:proofs:Equivalence}.
\end{proof}

{The above Lemma means that increasing the penalty coefficient makes the KKT solution of \eqref{eq:pen_reform} more feasible with respect to \eqref{eq:disc_reform}—that is, the term $\sum_{i=1}^{k} p_{i}^{\ast}g\left(\x^{\ast},\y^{(i)}\right) - \widetilde{V}_{\lambda}(\x^{\ast})$ becomes smaller. This result plays a key role in analyzing the equivalence between \eqref{eq:disc_reform} and \eqref{eq:pen_reform}, as it allows us to evaluate the range of relax coefficient values in \eqref{eq:disc_reform} for which the solution to \eqref{eq:pen_reform} remains feasible for \eqref{eq:disc_reform}. This is also a standard conclusion in penalty methods, e.g., in \cite[Proposition 4]{shen2023penaltybased} the authors proved a similar result.}

Then, we show that these KKT points are not only feasible, but they are in fact KKT points of \eqref{eq:disc_reform}.

\begin{proposition}\label{pro:disc_pen_optima}
Suppose Assumptions \ref{ass:basics}, \ref{ass:cover} hold and $\lambda \leq \frac{L_f D}{\gamma}$. If $(\x^{\ast},\p^{\ast})$ is a KKT point or local/global minimum of \eqref{eq:pen_reform}, then there exists $\epsilon\leq\frac{5L_f D}{2\gamma}$ such that $(\x^{\ast},\p^{\ast})$ is also a KKT point or local/global minimum, respectively, of \eqref{eq:disc_reform}.
\end{proposition}
\begin{proof}
    See Appendix \ref{app:proofs:Equivalence}.
\end{proof}
The above proposition demonstrates that if we obtain some point that ``solves" (i.e., a local/global minimum or a KKT point) the penalty problem \eqref{eq:pen_reform}, then this point is also a solution (in the corresponding sense) of the discretized reformulation \eqref{eq:disc_reform} that we ultimately aim to solve. This connection is crucial, especially regarding the correspondence between the KKT points. Since our first-order methods only guarantee convergence to a KKT point of \eqref{eq:pen_reform}, we must ensure that such a point is also meaningful for the original problem \eqref{eq:blo_vf}. In particular, it should at least be a KKT point of \eqref{eq:blo_vf}. In \cite[Proposition 3]{shen2023penaltybased}, the authors establish a relationship between the $\epsilon$-KKT point of the original bilevel problem and the KKT point of its penalized reformulation. However, their approach only applies to the case where the lower-level problem is unconstrained. In contrast, we conduct a detailed comparison of the KKT conditions for the two problems, and are able to establish a connection between their KKT points even when the lower-level problem is constrained. To the best of our knowledge, this is the first work in the bilevel setting that establishes a correspondence between the KKT points of a vlaued-function-based problem with constrained lower-level problem and its penalized reformulation.

\section{Algorithm and Convergence Results}\label{sec:alg}

Having established the equivalence among \eqref{eq:blo_vf}, \eqref{eq:disc_reform}, and \eqref{eq:pen_reform} in Section \ref{sec:EquivalenceBetweenDifferentProblems}, we proceed to solve \eqref{eq:pen_reform} using a standard projected gradient descent approach. We will show that it provably converges to a stationary solution. To this end, we first establish that the objective of the proposed penalty problem \eqref{eq:pen_reform} has Lipschitz continuous gradients. 
\begin{lemma}\label{lem:Lipschitz}
    Suppose Assumptions \ref{ass:basics}, \ref{ass:cover} hold. The objective function of the penalty reformulation \eqref{eq:pen_reform},
    $$F_{\gamma}(\x,\p)=f\left(\x,\sum_{i=1}^k p_i\y^{(i)}\right)+\gamma\left(\sum_{i=1}^k p_ig(\x,\y^{(i)})-\widetilde{V}_{\lambda}(\x)\right)$$
    is $L_\gamma$-Lipschitz-smooth {with respect to $(\x,\p)$} where 
    \begin{equation}
        L_\gamma=\left\{\overline{L}_f\max\{1,\sqrt{k}D^2\}+\gamma\left[\left(\left(\overline{L}_g+L_g\right)^2+kL_g^2\right)^{\frac{1}{2}}+\frac{L_g^2k}{\lambda}+\overline{L}_g\right]\right\}.
    \end{equation}
\end{lemma}
\begin{proof}
    {See} Appendix \ref{app:proofs:alg}.
\end{proof}
The Lipschitz smoothness of the objective of \eqref{eq:pen_reform} enables the development of a provably convergent {projected} gradient-descent-based method named DIscretized Value-FunctIon-BaseD - Bilevel Optimization (DIVIDE-BLO) {for efficiently solving the problem}. DIVIDE-BLO is an iterative algorithm over the variables $(\x,\p)$ that primarily performs the following two steps. First, it computes the value of the approximate value function $\widetilde{V}(\x)$ at the current iterate $\x$. Note that computing this value does not involve solving a generic strongly convex minimization problem, as in the case of the typical value function $V(\x)$. Instead, the minimization of the problem in $\widetilde{V}(\x)$ reduces to finding the projection of the point $\left(-\frac{g(\x,\y^{(1)})}{\lambda},\cdots,-\frac{g(\x,\y^{(k)})}{\lambda} \right)^\top $ {onto} the unit simplex; for the relevant explanation please check the derivations in \eqref{eq:vtilde_proj}, in the proof of Proposition \ref{pro:dvf_properties}. Next, DIVIDE-BLO performs one projected gradient descent step on $F_{\gamma}(\x,\p)$ over the pair of variables $(\x,\p)$, using the previously obtained value $\widetilde{V}(\x)$. For a complete description of the proposed DIVIDE-BLO algorithm, see Algorithm \ref{alg}.

\begin{algorithm}[t!]
\caption{\underline{Di}scretized \underline{V}alue Funct\underline{i}on Base\underline{d} M\underline{e}thod for \underline{B}i\underline{l}evel \underline{O}ptimization (DIVIDE-BLO)}
\begin{algorithmic}
\STATE \textbf{Step 0: Initialization.} 
\begin{itemize}
    \item[-] Sample $k$ points $\widetilde{\mathcal{Y}}=\{\y^{(i)}\}_{i=1}^{k} \in \mathcal{Y}$. 
    \item[-] Select the penalty parameter $\gamma$, the regularization constant $\lambda$, the accuracy level $\epsilon$, and the stepsize $\alpha$. 
    \item[-] Generate an initial point $\x_0$ and $\p_0$. Set $t = 0$.
\end{itemize}

\STATE \textbf{Step 1: Evaluate value function.} 
\begin{itemize} 
    \item[-]  Compute $g(\x_t,\y^{(i)})$ for each $i \in [k]$ and set $\hat{\p}=\left(-\frac{g(\x_t,\y^{(1)})}{\lambda},\ldots,-\frac{g(\x_t,\y^{(k)})}{\lambda}\right).$
    \item[-] Compute the projection $\p^{\ast}=\text{Proj}_{\Delta}(\hat{\p})$. 
    
    Set $\widetilde{V}_{\lambda}(\x_t)=\sum_{i=1}^k p^{\ast}_i g(\x_t,\y^{(i)})+\frac{\lambda}{2}\left\|\p^{\ast}\right\|^2$ 
\end{itemize}
\STATE \textbf{Step 2: Compute the gradient of $F_\gamma(\x_t,\p_t)$.} 
{
\begin{align*}
    \begin{bmatrix}
        \nabla_{x} F_\gamma (\x_t,\p_t) \\
        \nabla_{p} F_\gamma (\x_t,\p_t)
    \end{bmatrix}
    =
    \begin{bmatrix}
        \nabla_{x} \widetilde{f}(\x_t,\p_t) \\
        \nabla_{p} \widetilde{f}(\x_t,\p_t)
    \end{bmatrix}
    + \gamma
    \begin{bmatrix}
        \nabla_{x} \widetilde{g}(\x_t,\p_t) - \sum_{i=1}^kp_{i}^{\ast}\nabla_{\x} g(\x,\y^{(i)}) \\
        \nabla_{p} \widetilde{g}(\x_t,\p_t)
    \end{bmatrix}
\end{align*}%
}
where $\widetilde{f}(\x,\p)=f\left(\x,\sum_{i=1}^k p_i\y^{(i)}\right)$ and $\widetilde{g}(\x,\p)=\sum_{i=1}^k p_i g\left(\x,\y^{(i)}\right)$. 
\STATE \textbf{Step 3: Update $\x_{t+1}$ and $\p_{t+1}$.}
\begin{align*}
(\x_{t+1},\p_{t+1})=\text{Proj}_{\mathcal{X}\times\Delta}\left(\left(\x_t,\p_t\right)^\top-\alpha\nabla F_\gamma (\x_t,\p_t)\right) 
\vspace{-2 mm}
\end{align*}
\STATE \textbf{Step 4: Termination check.} 

If the following condition is satisfied, 
\begin{align*}
    \frac{1}{\alpha}\left\|(\x_{t},\p_{t})-(\x_{t+1},\p_{t+1})\right\|\leq\epsilon,
\end{align*}
then {\bf Stop}. Otherwise, let $t=t+1$, and return to {\bf Step 1}.
\end{algorithmic}
\label{alg}
\end{algorithm}

Next, we show the convergence of DIVIDE-BLO to a stationary solution. To this end, and taking into account the fact that the problem we are solving is constrained, we use the following projected gradient to {quantify} the convergence
of problem \eqref{eq:pen_reform}:
\begin{align*}
    \left\|\nabla \mathcal{G}_{t}^{\alpha}\right\| = \frac{1}{\alpha}\left\|(\x_t,\p_t)-\text{Proj}_{\mathcal{X}\times\Delta}\left((\x_t,\p_t)-\alpha\nabla F_\gamma(\x_t,\p_t)\right)\right\|.
\end{align*}
Below we present the convergence {guarantees achieved by} DIVIDE-BLO.
\begin{theorem}\label{thm:convergence}
    Suppose Assumptions \ref{ass:basics}, \ref{ass:cover} hold, and select $\alpha\leq L_\gamma^{-1}$. Then the following hold for the iterates {generated by the} DIVIDE-BLO (Algorithm \ref{alg}):
    \begin{enumerate}[leftmargin = *]
        \item $\frac{1}{T}\sum\limits_{i=1}^{T}\left\|\nabla \mathcal{G}_{t}^{\alpha}\right\|^2\leq 2\frac{F_\gamma(\x_1,\p_1)-C_f+\frac{\lambda}{2}}{\alpha T}$, where $C_f=\inf_{(\x,\p)\in\mathcal{X}\times\Delta}\widetilde{f}(\x,\p)$.
        \item If the iterates $\{(\x_t,\p_t)\}_{t=1}^T$ converge, then they converge to a KKT point of \eqref{eq:pen_reform}, and therefore (according to {Proposition} \ref{pro:disc_pen_optima}) to a KKT point of \eqref{eq:disc_reform} for some {$\epsilon\leq\frac{5L_f D}{2\gamma}$}. If the KKT point is a global/local optimum of \eqref{eq:pen_reform}, then it is also a global/local optimum, respectively, of \eqref{eq:disc_reform}.
    \end{enumerate}
\end{theorem}
\begin{proof}
    {See} Appendix \ref{app:proofs:alg}.
\end{proof}
The above theorem establishes {finite-time} convergence of DIVIDE-BLO {to} a stationary point of {the problem} \eqref{eq:pen_reform} at a {rate of} $\frac{1}{T}$. In addition, {we have established} that asymptotically the iterates {generated by} DIVIDE-BLO attain a KKT point of the penalty reformulated problem \eqref{eq:pen_reform}. This result, combined with the connections established between the problems \eqref{eq:pen_reform} and \eqref{eq:disc_reform} in {Proposition} \ref{pro:disc_pen_optima}, implies that the iterates of DIVIDE-BLO also asymptotically attain a KKT point of \eqref{eq:disc_reform}, for some suitable {choice of $\epsilon > 0$}.

\section{Experiments}

\subsection{Numerical Experiments}

We consider the following problem:
\begin{align}\label{numerical_example}
    &\min\limits_{x \in [-10,10]}   f(x,y^\ast)=3x^2+7xy^\ast+5(y^\ast)^2+x-y+5  \\
    & \text{s.t. } y^\ast \in \arg\min\limits_{y \in [-5,5]} g(x,y)=\sin(10y)+2y^2.\nonumber 
\end{align}
Note that the lower-level objective function $g(x,y)$ is highly non-convex in $y$. We compare the proposed DIVIDE-BLO algorithm with the following baselines: 1) TTSA algorithm \cite{hong2020two} (an implicit gradient-based method), 2) V-PBGD algorithm \cite{shen2023penaltybased} (a penalty-based method). We conducted experiments with the following settings:
\begin{itemize}[leftmargin = *]
    \item DIVIDE-BLO: The initial point $x_0$ is uniformly chosen from $[-10, 10]$. The number of sampling points is $K = 201$, evenly sampled in $[-5, 5]$, forming a covering with radius $0.05$. The initial weight vector $\mathbf{p}_0$ is randomly initialized on the simplex. The penalty parameter $\gamma$ is selected from $\{20, 40, 60, 80, 100\}$, and the step size $\alpha$ is selected from $\{0.1, 0.05, 0.01\}$.

    \item TTSA: The initial points $x_0$ and $y_0$ are uniformly chosen from $[-10, 10]$ and $[-5, 5]$, respectively. Time-decaying step sizes are used: $\alpha_1(t) = \eta_1 / t^{3/5}$ and $\alpha_2(t) = \eta_2 / t^{2/5}$, where $\eta_1, \eta_2 \in \{0.1, 0.05, 0.01\}$.

    \item V-PBGD: The initial points $x_0$ and $y_0$ are uniformly chosen from $[-10, 10]$ and $[-5, 5]$, respectively. The penalty parameter $\gamma$ is selected from $\{20, 40, 60, 80, 100\}$. The outer-loop step size $\alpha_1$ is chosen from $\{0.1, 0.05, 0.01\}$, and the inner-loop step size $\alpha_2$ is chosen from $\{0.01, 0.005, 0.001\}$.
\end{itemize}

We define the violation of lower-level problem as 
\begin{align*}
    g(x_T,y_T)-\min_{y}g(x_T,y), 
\end{align*}
where $(x_T,y_T)$ is the output of the DIVIDE-BLO, TTSA, V-PBGD algorithms at iteration $T$. We also denote the overall gap as 
\begin{align*}
    (f(x_T,y_T)-f^\ast)+(g(x_T,y_T)-\min_{y}g(x_T,y)),
\end{align*}
where $f^\ast$ is the optimal value of problem (\ref{numerical_example}). 

We repeat the optimization process $50$ times for each parameter combination. For \textbf{DIVIDE-BLO}, we observe that for every combination of step size $\alpha$ and penalty parameter $\gamma$, and any initial point, the algorithm consistently converges to a solution with zero lower-level violation and zero total gap. In contrast, for both \textbf{TTSA} and \textbf{V-PBGD}, we find that for every combination of parameters, there always exists at least one initialization that leads to a solution with nonzero violation of the lower-level problem. 

\begin{figure}[htbp]
    \centering
    \includegraphics[width=0.45\textwidth, trim=0 0 0 20, clip]{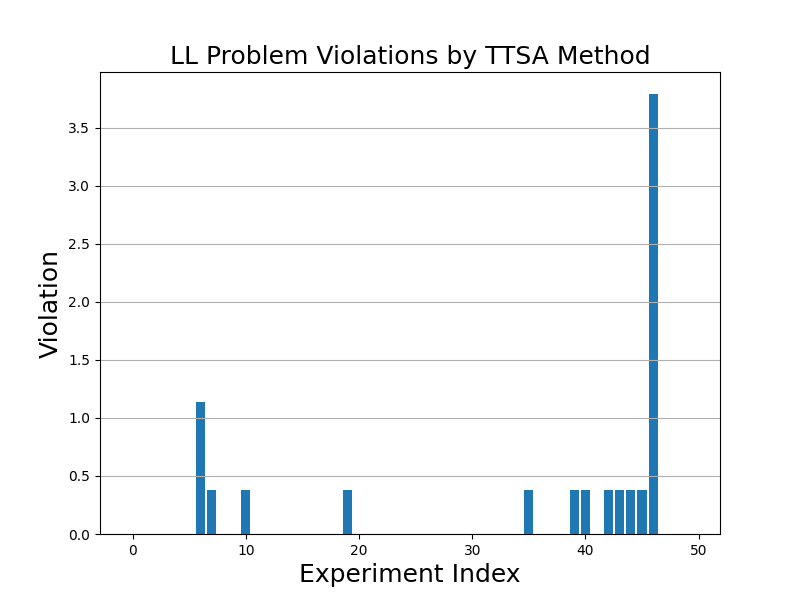}
    \includegraphics[width=0.45\textwidth, trim=0 0 0 20, clip]{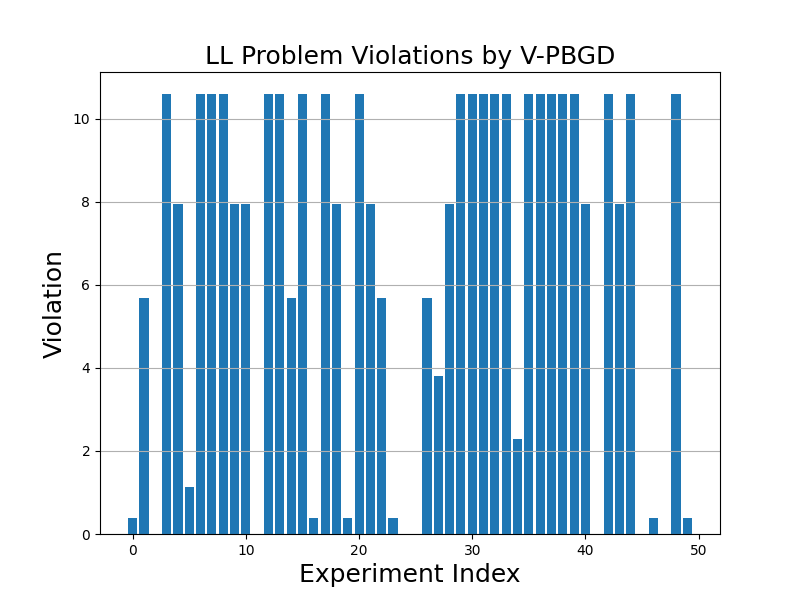}
    \includegraphics[width=0.45\textwidth, trim=20 0 0 20, clip]{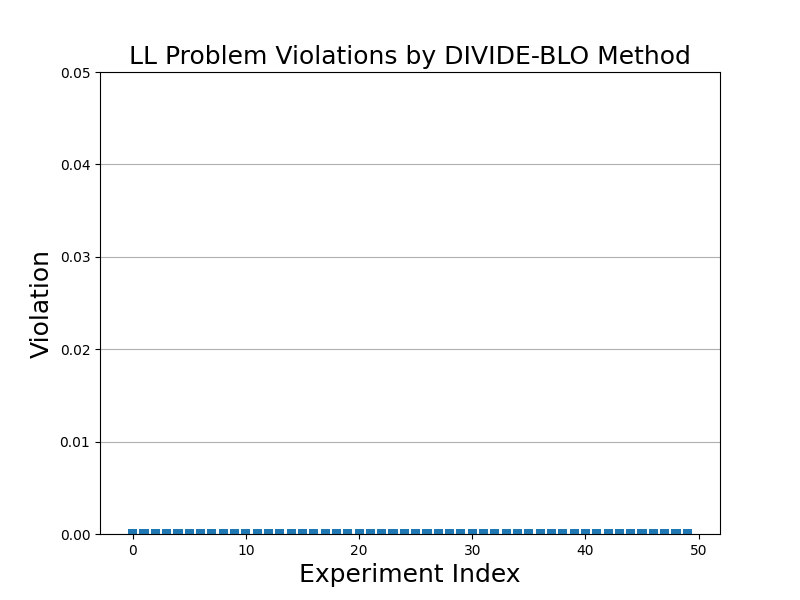}
    \includegraphics[width=0.42\textwidth, trim=0 0 0 10, clip]{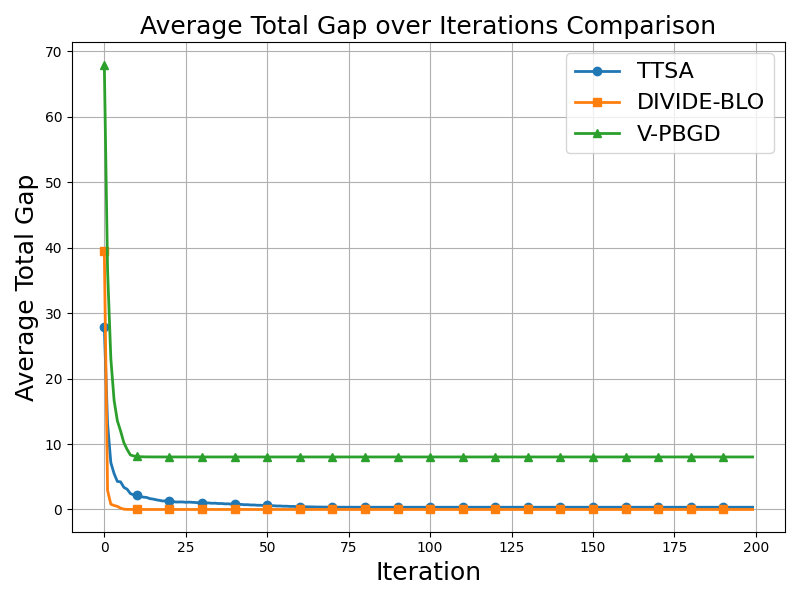}
    \caption{The top two figures and the bottom left figure respectively show the lower-level violations from 50 runs under the best parameter settings for TTSA, V-PBGD, and DIVIDE-BLO. The figure in the bottom right corner shows the average total gap curves over iterations for each algorithm, averaged across 50 randomly initialized runs. Only DIVIDE-BLO guarantees convergence to solutions with zero lower-level violation and zero total gap.
}
    \label{fig:num_exp}
\end{figure}

We selected the best parameter setting for each algorithm based on the lowest average total gap across all 50 runs. Figure~\ref{fig:num_exp} shows the lower-level violations under these optimal settings, confirming that DIVIDE-BLO consistently achieves zero violation across all runs, while TTSA and V-PBGD do not. We also plot the average total gap curves over iterations for each algorithm under its best configuration. The selected parameters are: DIVIDE-BLO with step size $\displaystyle \alpha = 0.1$ and penalty parameter $\displaystyle \gamma = 20$; TTSA with inner‐loop step size $\displaystyle \alpha_1 = 0.1/t^{3/5}$ and outer‐loop step size $\displaystyle \alpha_2 = 0.1/t^{2/5}$; V-PBGD  with penalty parameter $\displaystyle \gamma = 80$, inner‐loop step size $\displaystyle \alpha_1 = 0.05$, and outer‐loop step size $\displaystyle \alpha_2 = 0.005$.

The results indicate that the DIVIDE-BLO algorithm can effectively drive both the violation and total gap toward zero, ensuring that the attained solution of the lower-level problem is indeed close to the global optimum. TTSA and V-PBGD are not able to do so  because, due to the non-convexity of the lower-level problem, the  algorithms often get stuck at local optima, rendering the solutions infeasible with respect to the constraint $y^\ast \in \arg\min_{y \in [-5,5]} g(x,y)=\sin(10y)+2y^2$. This numerical experiment demonstrates the remarkable effectiveness of the DIVIDE-BLO algorithm in obtaining feasible solutions when the lower-level problem is non-convex.

\begin{table}[t]
\centering
\caption{The performance (classification accuracy) of the base and the weighted models across three datasets Cifar10, SVHN, and FashionMNIST, and three different pools of base models (Strong, Weak, and Weak+Strong). ``Base'' refers to the performance attained by the individual base models, and ``DIVIDE-BLO'' refers to the performance attained by the weighted model obtained with the use of the DIVIDE-BLO algorithm (for solving problem \eqref{eq:bilevel_ensemble}).}
\begin{tabular}{|c|c|c|c|c|}
\hline
Dataset & Model & Strong & Weak & Strong+Weak \\ \hline
\multirow{2}{*}{Cifar10} & Base & \makecell{$45.2$, $48.6$, $49.2$,\\ $50.4$, $50.9$, $51.3$, \\ $51.9$, $52.5$, $53.1$, \\  $53.3$ } &  \makecell{$20.0$, $20.3$, $23.0$, \\ $21.4$, $23.5$, $24.2$, \\ $26.8$, $27.1$, $27.2$, \\ $29.6$, $31.7$, $32.3$, \\ $32.5$, $33.2$}    &  \makecell{$24.2$, $26.8$, $27.2$, \\ $31.7$, $32.3$, $33.2$, \\ $48.6$, $49.2$, $51.9$, \\ $51.9$, $52.5$, $53.3$} \\ \cline{2-5} 
& DIVIDE-BLO & $\mathbf{59.2}$ & $\mathbf{41.0}$ & $\mathbf{57.3}$ \\ \hline
\multirow{2}{*}{SVHN} & Base & 
\makecell{$76.7$, $77.1$, $81.2$,\\ $82.9$, $83.5$, $84.04$,\\ $84.4$, $84.5$, $84.6$, \\ $84.9$, $85.1$, $85.2$, \\ $85.2$, $86.2$, $86.7$} & \makecell{$19.7$, $19.8$, $22.8$, \\ $28.5$, $29.2$, $38.4$, \\ $38.6$, $50.7$, $51.8$,\\ $57.1$} & \makecell{$28.5$, $38.4$, $38.6$,\\ $50.7$, $51.8$, $57.1$, \\ $81.2$, $84.4$, $84.9$,\\$85.1$, $85.2$, $86.7$} \\ \cline{2-5} 
& DIVIDE-BLO   & $\mathbf{90.0}$ & $\mathbf{65.8}$ & $\mathbf{89.3}$ \\ \hline
\multirow{2}{*}{\makecell{Fashion \\ MNIST}} & Base & \makecell{$75.0$, $75.2$, $77.4$, \\ $78.3$, $79.9$, $80.8$, \\ $81.8$, $82.2$, $82.6$, \\ $83.6$, $84.7$, $85.1$} & \makecell{$10.4$, $14.9$, $17.2$, \\ $19.2$, $20.3$, $28.0$, \\ $44.0$, $51.5$, $52.2$, \\ $58.3$} & \makecell{$14.9$, $17.2$, $20.3$, \\$44.0$, $51.5$, $52.2$, \\ $79.9$, $81.8$, $82.2$, \\$82.6$, $84.7$, $85.1$} \\ \cline{2-5} 
& DIVIDE-BLO &  $\mathbf{85.5}$ & $\mathbf{66.8}$ &  $\mathbf{85.8}$ \\ \hline
\end{tabular}
\label{tab:ensemble}
\end{table}

\subsection{Ensemble Learning}\label{sec:exp_ensemble}
We consider the ensemble learning problem introduced in Sec. \ref{sec:apps} where the underlying learning problem is a classification task. The goal is through the solution of problem \eqref{eq:bilevel_ensemble} to obtain a weighted model that performs better than the individual base models on the same task. The experimental procedure is as follows. {First, we generate multiple models by repeatedly pretraining a neural network using different training parameter configurations.} This produces a set of models with varying performance levels (classification accuracy), from weak to strong. Using this set we create three different pools of base models. A pool that only has strong models (high performance), a pool that only has weak models (low performance), and a pool that has both weak and strong models. 

For obtaining a high-performing weighted model we use the DIVIDE-BLO algorithm for solving the bilevel problem \eqref{eq:bilevel_ensemble}. We use $5000$ and $2000$ samples for $D^{train}$ and $D^{valid}$, respectively, while the test set on which all the models (individuals and weighted) are evaluated consists of $5000$ samples. We consider three datasets: CIFAR10, SVHN and FashionMNIST. The key parameter values, i.e., the number of points $k$ of the discrete set $\widetilde{\mathcal{Y}}$, the stepsize $\alpha$, the regularization parameter $\lambda$, the penalty parameter $\gamma$ and the number of iterations are selected with the use of a hyperparameter optimization procedure. For example,  the number of points $k$ is selected within the range $[10,40]$. In Table \ref{tab:ensemble} we report the performance of the base models and the performance of the weighted model (referred to as ``DIVIDE-BLO'' in the table) we obtain by solving the bilevel problem \eqref{eq:bilevel_ensemble}.
{We observe that the weighted model outperforms the base model in all experimental configurations.} This showcases the utility of the bilevel ensemble learning formulation and the ability of DIVIDE-BLO to find good solutions to such problems.

\section{Conclusion}
In this work we developed a value function and penalty-based method for the solution of non-convex bilevel problems where the LL has constraints. The key idea is the discretization of the value function that results to a reformulated problem that possesses certain attractive properties (e.g., smoothness). The DIVIDE-BLO algorithm is proposed for the solution of the reformulated problem and its convergence is established. Nonetheless, there are opportunities for further analysis. For instance, we will be interested in studying the effect and sample complexity of different, potentially adaptive, sampling mechanisms in the algorithm's performance. Another direction of interest to us is the identification of more applications with low-dimensional LL bilevel structure in areas such as machine learning, operations research and signal processing.

\appendix

\section{Proofs}

\subsection{Proofs of Section \ref{sec:reform}}\label{app:proofs:reform}

\begin{proof}[\textbf{Proof of Proposition \ref{pro:dvf_properties}}]

\hfill \break
\textbf{(a)} The value function $\widetilde{V}_{\lambda}(\x)$ can be rewritten as
\begin{align*}
    \widetilde{V}_{\lambda}(\x)&=\min\limits_{{\bf 0} \leq \p \leq {\bf 1}, {\bf 1}^\top\p=1} \left\{ \sum\limits_{i=1}^{k} p_{i}g(\x,\y^{(i)}) + \frac{\lambda}{2}\|\p\|^{2}\right\} \\
    &=-\max\limits_{{\bf 0} \leq \p \leq {\bf 1}, {\bf 1}^\top\p=1} \left\{ -\sum\limits_{i=1}^{k} p_{i}g(\x,\y^{(i)}) - \frac{\lambda}{2}\|\p\|^{2}\right\}.
\end{align*}
Note that the constraint set of the above maximization problem is the unit simplex, which is a convex set. Furthermore, the objective is strongly concave, and as a result the respective maximization problem has a unique solution. It then follows from Danskin's theorem \cite[Corollary, pg. 1167]{bernhard1995danskin} that the approximate value function $\widetilde{V}_{\lambda}(\x)$ is differentiable and its gradient is given by the formula,
\begin{align*}
    \nabla \widetilde{V}_{\lambda}(\x) = \sum\limits_{i=1}^{k} p_{i}^{\ast}\nabla_{\x}g(\x,\y^{(i)}) 
\end{align*}
where $p_{i}^{\ast}$ are the elements of the optimal vector $$\p^{\ast} = \arg\max\limits_{{\bf 0} \leq \p \leq {\bf 1}, {\bf 1}^\top\p=1} \left\{ -\sum\limits_{i=1}^{k} p_{i}g(\x,\y^{(i)}) - \frac{\lambda}{2}\|\p\|^{2}\right\}.$$

\textbf{(b)} The following hold: 
\begin{align*}
\left\lvert \widetilde{V}_{\lambda}(\x)-V(\x)\right\rvert
&\stackrel{(a)}{=}\min\limits_{{\bf 0} \leq \p \leq {\bf 1}, {\bf 1}^\top\p=1} \left\{ \sum\limits_{i=1}^{k} p_{i}g(\x,\y^{(i)})+ \frac{\lambda}{2}\|\p\|^{2}\right\}-\min\limits_{\y \in \mathcal{Y}} \left\{ g(\x,\y)\right\}\\
&\stackrel{(b)}{\leq} g(\x,\y^{(j)})- g(\x,\y^{\ast}) + \frac{\lambda}{2}\\
&\stackrel{(c)}{\leq} L_{g} \|\y^{(j)}-\y^{\ast}\| + \frac{\lambda}{2}\\
&\stackrel{(d)}{\leq} \frac{2L_gD\sqrt{m}}{k^{1/m}}+\frac{\lambda}{2}.
\end{align*}
The equality in (a) follows from the fact that $\widetilde{V}_{\lambda}(\x)\geq V(\x)$; in (b) we denote with $\y^{\ast} \in \mathcal{Y}$ one of the global minima of $V(\x)$ and we set $p_{j}=1,p_{i}=0, \text{ for } i \neq j$, where  $j=\arg\min_{i=1,\ldots,k} \|\y^{(i)}-\y^{\ast}\|$; in (c) we make use of the Lipschitz continuous property of $g$ (Assumption \ref{ass:basics}\ref{ass:basics:Lip_bound_g}). Finally, in (d) we know that for a given bounded set $\mathcal{Y} \in \mathbb{R}^{m}$ with $\|\y\| \leq D, \forall \y \in \mathcal{Y}$ and a covering of radius $r>0$ it holds that \cite{shalev2014understanding} $$k<\left(\frac{2D\sqrt{m}}{r}\right)^{m} \implies r<\frac{2D\sqrt{m}}{k^{1/m}}.$$ 

\textbf{(c)}
To begin with, we note that the problem in \eqref{eq:DVF} can be rewritten as
\begin{align}\label{eq:vtilde_proj}
    \p^{\ast}(\x)&:=\arg\min\limits_{{\bf 0} \leq \p \leq {\bf 1}, {\bf 1}^\top\p=1} \left\{ \sum\limits_{i=1}^{k} p_{i}g(\x,\y^{(i)}) + \frac{\lambda}{2}\|\p\|^{2}\right\} \nonumber\\
    &=\arg\min\limits_{{\bf 0} \leq \p \leq {\bf 1}, {\bf 1}^\top\p=1} \left\{ \sum\limits_{i=1}^{k} \left(p_{i}g(\x,\y^{(i)}) + \frac{\lambda}{2}p_i^{2}\right)\right\} \nonumber\\
    &=\arg\min\limits_{{\bf 0} \leq \p \leq {\bf 1}, {\bf 1}^\top\p=1} \left\{\sum_{i=1}^k\left(p_i+\frac{g(\x,\y^{(i)})}{\lambda}\right)^2\right\} \nonumber\\
    &=\arg\min\limits_{{\bf 0} \leq \p \leq {\bf 1}, {\bf 1}^\top\p=1} \left\{\left\|\p-\left[-\frac{g(\x,\y^{(1)})}{\lambda},\cdots,-\frac{g(\x,\y^{(k)})}{\lambda} \right]^\top  \right\|^{2}\right\}  \nonumber\\
    &=\text{proj}_{{\bf 0} \leq \p \leq {\bf 1}, {\bf 1}^\top \p=1}\left(\left[-\frac{g(\x,\y^{(1)})}{\lambda},\cdots,-\frac{g(\x,\y^{(k)})}{\lambda} \right]^\top \right),
\end{align}
where in the last expression we see that $\p^{\ast}(\x)$ is the projection  of $$\hat{\p}(\x):=\left[-\frac{g(\x,\y^{(1)})}{\lambda},\cdots,-\frac{g(\x,\y^{(k)})}{\lambda} \right]^\top $$ onto the unit simplex.  

Then, by the non-expansiveness of the projection operator, we have
\begin{align}\label{eq:proj_nonexpan}
    \|\p^{\ast}(\x_1)-\p^{\ast}(\x_2)\|&\leq\|\hat{\p}(\x_1)-\hat{\p}(\x_2)\| \nonumber\\
    &=\frac{1}{\lambda}\left\{\sum_{i=1}^k\left(g(\x_1,\y^{(i)})-g(\x_2,\y^{(i)})\right)^2\right\}^{\frac{1}{2}} \nonumber\\
    &\leq \frac{1}{\lambda}\left\{kL_g^2\|\x_1-\x_2\|^2\right\}^{\frac{1}{2}} \nonumber\\
    &=\frac{L_g\sqrt{k}}{\lambda}\|\x_1-\x_2\|,
\end{align}
where in the second inequality above we used the Lipschitz continuity of $g$ (Ass. \ref{ass:basics}\ref{ass:basics:Lip_bound_g}).

Therefore, we obtain
\begin{align*}
    \|\nabla\widetilde{V}_{\lambda}(\x_1)-\nabla\widetilde{V}_{\lambda}(\x_2)\|=&\left\|\sum\limits_{i=1}^{k} p_{i}^{\ast}(\x_1)\nabla_{\x}g(\x_1,\y^{(i)})-\sum\limits_{i=1}^{k} p_{i}^{\ast}(\x_2)\nabla_{\x}g(\x_2,\y^{(i)})\right\|\\
    \leq&\left\|\sum\limits_{i=1}^{k} p_{i}^{\ast}(\x_1)\nabla_{\x}g(\x_1,\y^{(i)})-\sum\limits_{i=1}^{k} p_{i}^{\ast}(\x_2)\nabla_{\x}g(\x_1,\y^{(i)})\right\|\\
    &+\left\|\sum\limits_{i=1}^{k} p_{i}^{\ast}(\x_2)\nabla_{\x}g(\x_1,\y^{(i)})-\sum\limits_{i=1}^{k} p_{i}^{\ast}(\x_2)\nabla_{\x}g(\x_2,\y^{(i)})\right\|\\
    \leq&\sum\limits_{i=1}^{k}\left| p_{i}^{\ast}(\x_1)-p_{i}^{\ast}(\x_2)\right|\left\| \nabla_{\x}g(\x_1,\y^{(i)})\right\|\\
    &+\sum\limits_{i=1}^{k}\left| p_{i}^{\ast}(\x_2)\right|\left\| \nabla_{\x}g(\x_1,\y^{(i)})-\nabla_{\x}g(\x_2,\y^{(i)})\right\|\\
    \leq&L_g\sqrt{k}\|\p^{\ast}(\x_1)-\p^{\ast}(\x_2)\|+\overline{L}_g\|\x_1-\x_2\|\\
    \leq&\left(\frac{L_g^2k}{\lambda}+\overline{L}_g\right)\|\x_1-\x_2\|
\end{align*}
Note that in the first equality, we applied the gradient formula we derived in the (a) part of this proposition. Also, we obtain the final expression with the application of expression \eqref{eq:proj_nonexpan} and the Lipschitz continuity and gradient property of $g$. The proof is completed.
\end{proof}
\subsection{Proofs of Section \ref{sec:EquivalenceBetweenDifferentProblems}}\label{app:proofs:Equivalence}
\begin{proof}[\textbf{Proof of Lemma \ref{lem:ll_bound}}]
    According to Proposition \ref{pro:dvf_properties}, we have
    $$\left|\widetilde{V}_\lambda(\x)-V(\x)\right|\leq\frac{2L_gD\sqrt{m}}{k^{1/m}}+\frac{\lambda}{2}.$$
    By the assumption $k\geq(4L_gD\sqrt{m}/\epsilon_1)^m$ and $\lambda\leq\epsilon_1$, it follows that
    \begin{align}\label{eq:lessthen2epsilon}
        \sum_{i=1}^k p_ig(\x,\y^{(i)})-V(\x)=\sum_{i=1}^k p_ig(\x,\y^{(i)})-\widetilde{V}_\lambda(\x)+\widetilde{V}_\lambda(\x)-V(\x)\leq 2\epsilon_1.
    \end{align}
    Set $I:=\{i:g(\x,\y^{(i)})-V(\x)\leq\sqrt{\epsilon_1}\}$. Since $k\geq(4L_gD\sqrt{m}/\epsilon_1)^m$, the set of $k$ sample points forms an cover of $\mathcal{Y}$ with radius $\epsilon_1/(2L_g)$, {so for any $\y^\ast\in \arg\min g(\x,\y)$, we can find a sample point, denoted as $\y^{(i)}$, such that $\|\y^{(i)}-\y^\ast\|\leq\epsilon_1/(2L_g)$. This implies 
    $$\min_{1\leq i\leq k}g(\x,\y^{(i)})-V(\x)=g(\x,\y^{(1)})-V(\x)\leq\epsilon_1/2\leq\sqrt{\epsilon_1}$$since $\epsilon_1\leq 1$.} Then the following holds
    \begin{align}
        g(\x,\sum_{i=1}^kp_i\y^{(i)})-V(\x)&=g(\x,\sum_{i\in I}p_i\y^{(i)}+\sum_{i\notin I}p_i\y^{(1)}+\sum_{i\notin I}p_i(\y^{(i)}-\y^{(1)}))-V(\x)\nonumber\\
        &\overset{(\text{i})}{\leq} g(\x,\sum_{i\in I}p_i\y^{(i)}+\sum_{i\notin I}p_i\y^{(1)})-V(\x)+2DL_g\sum_{i\notin I}p_i\nonumber\\
        &\overset{(\text{ii})}{\leq}\sqrt{\epsilon_1}+2DL_g\sum_{i\notin I}p_i.\label{ineq:g-V}
    \end{align}
    {In (i), we use Assumption \ref{ass:basics}\ref{ass:basics:Lip_bound_g}, i.e., the Lipschitz continuity assumption, and Assumption \ref{ass:basics}\ref{ass:basics:sets}, which states that the distance between any two points in the constraint set $\mathcal{Y}$ is at most $2D$. In (ii), we use the fact that $g(\x,\y^{(i)}) - V(\x) \leq \sqrt{\epsilon_1}$ for all $i \in I$, by the construction of the index set $I = \{i : g(\x,\y^{(i)}) - V(\x) \leq \sqrt{\epsilon_1} \}$. We also use that the level set $\{\y : g(\x,\y) - V(\x) \leq \sqrt{\epsilon_1} \leq \epsilon_0\}$ is convex, as guaranteed by Assumption~\ref{ass:reform}(\ref{ass:reform(1)}). Note that for any $i\notin I$, we have $g(\x,\y^{(i)})-V(\x)>\sqrt{\epsilon}$. The following holds} 
    $$\sum_{i\notin I}p_i\sqrt{\epsilon_1}\leq\sum_{i\notin I}p_i(g(\x,\y^{(i)})-V(\x))\overset{(\text{iii})}{\leq} 2\epsilon_1,$$
    {where we use \eqref{eq:lessthen2epsilon} in (iii).} Therefore, $\sum_{i\notin I}p_i\leq 2\sqrt{\epsilon_1}$. Plug this into (\ref{ineq:g-V}), we obtain
    $$g\left(\x,\sum_{i=1}^kp_i\y^{(i)}\right)-V(\x)\leq(1+4DL_g)\sqrt{\epsilon_1}.$$
\end{proof}

\begin{proof}[\textbf{Proof of Theorem \ref{thm:global}}]
     \textbf{Part 1:} Suppose $(\widetilde{\x},\widetilde{\p})$ is a global solution of (\ref{eq:disc_reform}) with relaxation coefficient $\epsilon_1$. By Lemma \ref{lem:ll_bound}, we know that $g(\widetilde{\x},\widetilde{\y})-V(\widetilde{\x})\leq(1+4DL_g)\sqrt{\epsilon_1}$. {Suppose} there exists another solution $(\x',\y')$ such that $g(\x',\y')-V(\x')\leq(1+4DL_g)\sqrt{\epsilon_1}$ and $f(\x',\y')< f(\widetilde{\x},\widetilde{\y})-\delta$. By the assumption that $\epsilon_1\leq[\epsilon_0/(1+4DL_g)]^2$, we have $g(\x',\y')-V(\x')\leq\epsilon_0$. Therefore, according to Assumption \ref{ass:reform}(\ref{ass:reform(2)}), we can find a point $\y^\ast\in\arg\min\limits_{\y\in\mathcal{Y}}g(\x',\y)$ such that
    \begin{align*}
        \|\y^\ast-\y'\|{\leq C\left|g(\x',\y')-V(\x')\right|^\beta}\leq C\left[\left(1+4DL_g\right)\sqrt{\epsilon_1}\right]^\beta.
    \end{align*}
    By the assumption that $k\geq(4L_gD\sqrt{m}/\epsilon_1)^m$, there exists a sample point $\y^{(i)}$ such that $\|\y^{(i)}-\y^\ast\|\leq\epsilon_1/(2L_g)$. This implies
    \begin{align*}
        \|\y^{(i)}-\y'\|\leq  \|\y^{(i)}-\y^\ast\|+ \|\y^\ast-\y'\|\leq C\left[\left(1+4DL_g\right)\sqrt{\epsilon_1}\right]^\beta+\frac{\epsilon_1}{2L_g}.
    \end{align*}
    By the Lipschitz continuity of $f(x)$, we have 
    \begin{align*}
        |f(\x',\y^{(i)})-f(\x',\y')|\leq L_f\cdot\left\{C\left[\left(1+4DL_g\right)\sqrt{\epsilon_1}\right]^\beta+\frac{\epsilon_1}{2L_g}\right\}=\delta,
    \end{align*}
    {which implies
    \begin{align*}
        f(\x',\y^{(i)})\leq f(\x',\y')+\delta.
    \end{align*}}
    Suppose $\p'$ is a $k$-dimensional vector with its $i$-th component equal to $1$ and all the other components equal to $0$. We obtain $\widetilde{f}(\x',\p')=f(\x',\y^{(i)})\leq f(\x',\y')+\delta< f(\widetilde{\x},\widetilde{\y})=\widetilde{f}(\widetilde{\x},\widetilde{\p})$. Note that 
    \begin{align*}
        \sum_{i=1}^k p_ig(\x',\y^{(i)})-\widetilde{V}_\lambda(\x')&=g(\x',\y^{(i)})-\widetilde{V}_\lambda(\x')\\
        &\leq|g(\x',\y^{(i)})-\widetilde{V}(\x')|+|\widetilde{V}(\x')-\widetilde{V}_\lambda(\x')|\\
        &\leq|g(\x',\y^{(i)})-g(\x',\y^\ast)|+\frac{\lambda}{2}\\
        &\leq \frac{L_g\epsilon_1}{2L_g}+\frac{\epsilon_1}{2}=\epsilon_1,
    \end{align*}
    which means $(\x',\p')$ is feasible. This contradicts the assumption that $(\widetilde{\x},\widetilde{\p})$ is a global solution of (\ref{eq:disc_reform}) with relaxation coefficient $\epsilon_1$. {Therefore, we conclude that $(\x,\y)$ is a $\delta$-global solution of \eqref{eq:blo_vf} with relaxation coefficient $(1+4DL_g)\sqrt{\epsilon_1}$.}\\
    
    \noindent\textbf{Part 2:} Suppose $(\widetilde{\x},\widetilde{\y})$ is a global solution of (\ref{eq:blo_vf}) with relaxation coefficient $\epsilon_1$. By the assumption that $k\geq(4L_gD\sqrt{m}/\epsilon_1)^m$, we can find a sample point $\y^{(i)}$ such that $\|(\widetilde{\x},\widetilde{\y})-(\widetilde{\x},\y^{(i)})\|\leq\epsilon_1/(2L_g)$. This also implies 
    \begin{align*}
        \left|g(\widetilde{\x},\y^{(i)})-\widetilde{V}_\lambda(\x)\right|&\leq\left|g(\widetilde{\x},\y^{(i)})-g(\widetilde{\x},\widetilde{\y})\right|+\left|g(\widetilde{\x},\widetilde{\y})-V(\widetilde{\x})\right|+\left|V(\widetilde{\x})-\widetilde{V}_\lambda(\widetilde{\x})\right|\\
        &\leq\frac{\epsilon_1}{2}+\epsilon_1+\frac{\epsilon_1}{2}+\frac{\epsilon_1}{2}\\
        &=\frac{5}{2}\epsilon_1,
    \end{align*}
    and $f(\widetilde{\x},\y^{(i)})\leq f(\widetilde{\x},\widetilde{\y})+L_f\epsilon_1/(2L_g)$. Define $\widetilde{\p}$ be a $k$-dimensional vector with its $i$-th component equal to $1$ and all the other components equal to $0$. {We have
    \begin{align}\label{eq:widetildef(x,p)}
        \widetilde{f}(\widetilde{\x},\widetilde{\p})=f(\widetilde{\x},\y^{(i)})\leq f(\widetilde{\x},\widetilde{\y})+\frac{L_f\epsilon_1}{2L_g}.
    \end{align}}
    We aim to show $(\widetilde{\x},\widetilde{\p})$ is a $\delta$-global solution of (\ref{eq:disc_reform}) with {relaxation} coefficient $5\epsilon_1/2$. {Suppose} there exists another solution $(\x',\p')$ such that $\widetilde{f}(\x',\p')<\widetilde{f}(\widetilde{\x},\widetilde{\p})-\delta$ and $\sum_{i=1}^kp'_ig(\x,\y^{(i)})-\widetilde{V}_\lambda(\x)\leq 5\epsilon_1/2$. According to Lemma \ref{lem:ll_bound}, $(\x',\y'):=(\x',\sum_{i=1}^kp'_i\y^{(i)})$ satisfies
    $$g(\x',\y')-V(\x')\leq\frac{5}{2}\left(1+4DL_g\right)\sqrt{\epsilon_1}.$$
    By the assumption that $\epsilon_1\leq[2\epsilon_0/(5+20DL_g)]^2$, we have 
    $$g(\x',\y')-V(\x')\leq\frac{5}{2}\left(1+4DL_g\right)\sqrt{\epsilon_1}\leq\epsilon_0.$$
    In view of Assumption \ref{ass:reform}(\ref{ass:reform(2)}), there exists $\y^\ast\in\arg\min_{\y\in\mathcal{Y}}g(\x',\y)$ such that $$\|\y^\ast-\y'\|{\leq C\left|g(\x',\y')-V(\x')\right|^\beta}\leq C\left[\frac{5}{2}\left(1+4DL_g\right)\epsilon_1\right]^\beta.$$
    By continuity of $g(\x,\y)$, we can find a point $\y''$ on the line segment connecting $\y'$ and $\y^\ast$ such that
    \begin{align}\label{eq:feasibility}
        g(\x',\y'')-V(\x')\leq\epsilon_1.
    \end{align}
    Since $\y''$ is between $\y'$ and $\y^\ast$, we also have 
    \begin{align}\label{eq:y''-y'}
        \|\y''-\y'\|\leq \|\y^\ast-\y'\|\leq C\left[\frac{5}{2}\left(1+4DL_g\right)\epsilon_1\right]^\beta.
    \end{align}
    However,
    \begin{align*}
        f(\x',\y'')&\overset{\text{(i)}}{\leq} f(\x',\y')+CL_f\left[\frac{5}{2}\left(1+4DL_g\right)\sqrt{\epsilon_1}\right]^\beta\\
        &=\widetilde{f}(\x',\p')+CL_f\left[\frac{5}{2}\left(1+4DL_g\right)\sqrt{\epsilon_1}\right]^\beta\\
        &\overset{\text{(ii)}}{<}\widetilde{f}(\widetilde{\x},\widetilde{\p})+CL_f\left[\frac{5}{2}\left(1+4D\L_g\right)\sqrt{\epsilon_1}\right]^\beta-\delta\\
        &\overset{\text{(iii)}}{\leq} f(\widetilde{\x},\widetilde{\y})+L_f\cdot\left\{C\left[\frac{5}{2}\left(1+4DL_g\right)\sqrt{\epsilon_1}\right]^\beta+\frac{\epsilon_1}{2L_g}\right\}-\delta\\
        &= f(\widetilde{\x},\widetilde{\y}).
    \end{align*}
    {In (i), we apply Assumption~\ref{ass:basics}\ref{ass:basics:Lip_grad_f}, which states the Lipschitz smoothness of $f(\x,\y)$, together with \eqref{eq:y''-y'}. In (ii), we utilize the assumption for contradiction that $\widetilde{f}(\x',\p') < \widetilde{f}(\widetilde{\x}, \widetilde{\p}) - \delta$. Step (iii) follows directly from \eqref{eq:widetildef(x,p)}.}
    
    By (\ref{eq:feasibility}), $(\x',\y'')$ is a feasible solution. This contradicts the fact that $(\widetilde{\x},\widetilde{\y})$ is a global solution of (\ref{eq:BLO}) with relaxation coefficient $\epsilon_1$. {Therefore, we conclude that $(\widetilde{\x},\widetilde{\p})$ is a $\delta$-global solution of \eqref{eq:disc_reform} with relaxation coefficient $5\epsilon_1/2$.}
\end{proof}

\begin{proof}[\textbf{Proof of Lemma \ref{lem:KKT1}}]
\hfill \break
Let us set $\widetilde{f}(\x,\p)=f\left(\x,\sum\limits_{i=1}^k p_i\y^{(i)}\right)$ and $G(\x,\p)=\sum\limits_{i=1}^{k} p_{i}g\left(\x,\y^{(i)}\right) - \widetilde{V}_{\lambda}(\x)$.
The KKT points of \eqref{eq:pen_reform} satisfy the following conditions:
{\small
\begin{align}\label{KKT:PEN}
\begin{pmatrix}
\nabla_\x \widetilde{f}(\x,\p) \\
\nabla_{p_1} \widetilde{f}(\x,\p)\\
\vdots\\
\nabla_{p_i} \widetilde{f}(\x,\p)\\
\vdots\\
\nabla_{p_k} \widetilde{f}(\x,\p)
\end{pmatrix}
+
\gamma \begin{pmatrix}
\nabla_\x G(\x,\p) \\
g(\x,\y^{(1)})\\
\vdots\\
g(\x,\y^{(i)})\\
\vdots\\
g(\x,\y^{(k)})
\end{pmatrix}
-
\begin{pmatrix}
\mathbf{0} \\
\mu^0_1\\
\vdots\\
\mu^0_i\\
\vdots\\
\mu^0_k
\end{pmatrix}
\begin{matrix}
\ \\
\ \\
\ \\
% \leftarrow i\text{-th}\\
\ \\
\ 
\end{matrix}+\begin{pmatrix}
\mathbf{0} \\
\mu^1_1\\
\vdots\\
\mu^1_i\\
\vdots\\
\mu^1_k
\end{pmatrix}
\begin{matrix}
\ \\
\ \\
\ \\
% \leftarrow i\text{-th}\\
\ \\
\ 
\end{matrix}
+\begin{pmatrix}
\mathbf{0} \\
\mu\\
\vdots\\
\mu\\
\vdots\\
\mu
\end{pmatrix}=0,
\end{align}
}%

\noindent
where $\mu_i^0 \geq 0$, $\mu_i^1 \geq 0$, $i=1,\ldots,k$, and $\mu$ are the Lagrange multipliers of the constraints ${\bf 0} \leq \p $, $\p \leq {\bf 1}$, and ${\bf 1}^\top \p=1$, respectively. 
Without loss of generality, suppose there exists some integer $l$ with $1\leq l\leq k$ such that $p_1,\ldots,p_l$ are non-zero and $p_{l+1},\ldots,p_k$ are zero. It then follows {from complementary slackness} that
\begin{align}\label{eq:lagrange_mult}
    &\mu_1^0=\ldots=\mu_l^0=0 \nonumber\\
    &\mu_2^1=\ldots=\mu_{k}^1=0 \text{ if } \ell=1 \text{ or } \mu_1^1=\ldots=\mu_{k}^1=0 \text{ if } \ell>1.
\end{align}
We define
\begin{align}\label{eq:G}
    \hat{G}(\x,\p):= \sum_{i=1}^{l}p_i\cdot\left(g(\x,\y^{(i)})-g(\x,\y^{(s^*)})\right).
\end{align} 
where $s^*=\arg\min_s g(\x,\y^{(s)})$.

Then, depending on the exact value of the index $s^*$, we consider two different cases.

\noindent\textbf{Case 1:} $s^*\in\{l+1,..,k\}$
    
    By \eqref{KKT:PEN} we have
    \begin{align}
        \frac{\partial \widetilde{f}(\x,\p)}{\partial p_i}+\gamma g(\x,\y^{(i)})-\mu_i^0+\mu_i^1+\mu=0&\ \ \ \ \forall i\in\{1,...,l\}\label{eq:kkt1}\\
        \frac{\partial \widetilde{f}(\x,\p)}{\partial p_j}+\gamma g(\x,\y^{(j)})-\mu_j^0+\mu_j^1+\mu=0&\ \ \ \ \forall j\in\{l+1,...,k\}\label{eq:kkt2}
    \end{align}
    Subtracting eq. \eqref{eq:kkt2} from \eqref{eq:kkt1}, and using the conditions in \eqref{eq:lagrange_mult} and the nonnegativity of the Lagrange multipliers $\mu_i^0,\mu_i^1 \geq 0, i=1,\ldots,k$, we obtain the following,
    \begin{align}\label{eq:kkt_bound0}
        \gamma \left(g(\x,\y^{(i)})-g(\x,\y^{(j)})\right)\leq\frac{\partial \widetilde{f}(\x,\p)}{\partial p_j}-\frac{\partial \widetilde{f}(\x,\p)}{\partial p_i}.
    \end{align}
    We note that 
    \begin{align}\label{eq:trivialresult}
        \left|\frac{\partial \widetilde{f}(\x,\p)}{\partial p_j}\right|=\left|\nabla_{y}^\top f(\x,\y) \y^{(j)}\right| \leq \|\nabla_{y} f(\x,\y)\| \|\y^{(j)} \| \leq L_f D, j\in\{1,...,k\} 
    \end{align}
    due to assumptions \ref{ass:basics}\ref{ass:basics:sets}, \ref{ass:basics}\ref{ass:basics:Lip_bound_f} and Cauchy-Schwarz inequality. Then, applying this bound in \eqref{eq:kkt_bound0} we get 
    \begin{align*}
        g(\x,\y^{(i)})-g(\x,\y^{(j)})\leq\frac{2L_f D}{\gamma}.
    \end{align*}
    Especially, by setting $j=s^{\ast}$ in the above inequality, we have
    \begin{align}\label{eq:kkt_bound1}
        g(\x,\y^{(i)})-g(\x,\y^{(s^*)})\leq\frac{2L_f D}{\gamma}, \;\; i\in\{1,...,l\}.
    \end{align}
    Therefore, using \eqref{eq:G} we get
    \begin{align}\label{eq:G2_bound}
    \hat{G}(\x,\p)&=\sum\limits_{i=1}^lp_i\left(g\left(\x,\y^{(i)})-g(\x,\y^{(s^*)}\right)\right)\nonumber\\
    &\leq\sum\limits_{i=1}^lp_i\frac{2L_f D}{\gamma} \nonumber\\
    &=\frac{2L_f D}{\gamma},
    \end{align}
    where the inequality follows from \eqref{eq:kkt_bound1}. Note that
    \begin{align*}
        g(\x,\y^{(s^{\ast})})=\min\limits_{{\bf 0} \leq \p \leq {\bf 1}, {\bf 1}^\top \p=1} \sum_{i=1}^k p_ig(\x,\y^{(i)}),
    \end{align*}
    thus
    \begin{align*}
        \widetilde{V}_{\lambda}(\x)&=\min\limits_{{\bf 0} \leq \p \leq {\bf 1}, {\bf 1}^\top \p=1} \left\{ \sum\limits_{i=1}^{k} p_{i}g\left(\x,\y^{(i)}\right) + \frac{\lambda}{2}\|\p\|^{2}\right\}\\
        &\geq \min\limits_{{\bf 0} \leq \p \leq {\bf 1}, {\bf 1}^\top \p=1} \sum_{i=1}^k p_ig(\x,\y^{(i)})-\max\limits_{{\bf 0} \leq \p \leq {\bf 1}, {\bf 1}^\top \p=1} \sum_{i=1}^k \frac{\lambda}{2}\|\p\|^{2}\\
        &=g(\x,\y^{(s^{\ast})})-\frac{\lambda}{2}.
        \end{align*}
    
    Finally, with the use of the above inequality, we obtain the following bound
    \begin{align*}
        \sum_{i=1}^kp_ig(\x,\y^{(i)})-\widetilde{V}_{\lambda}(\x)&=\sum_{i=1}^kp_i\left(g(\x,\y^{(i)})-g(\x,\y^{(s^{\ast})})\right)\\
        &\qquad +g(\x,\y^{(s^{\ast})})-\widetilde{V}_{\lambda}(\x)\\
        &\qquad \leq \frac{2L_f D}{\gamma}+\frac{\lambda}{2}\\
        &\qquad =\frac{5L_fD}{2\gamma}.
    \end{align*}
    
\noindent\textbf{Case 2:} $s^*\in\{1,...,l\}$

    We subtract eq. \eqref{eq:kkt1} with $i=s^{\ast}$ from eq. \eqref{eq:kkt1} with $i \in \{1,\ldots,l\}$
    \begin{align*}
        &g(\x,\y^{(i)})-g(\x,\y^{(s^*)})\nonumber\\
       =&\frac{\mu_i^0-\mu_i^1-\mu_{s^*}^{0}+\mu_{s^*}^{1} 
         +\left(\frac{\partial \widetilde{f}(\x,\p)}{\partial p_{s^*}}-\frac{\partial \widetilde{f}(\x,\p)}{\partial p_i}\right)}{\gamma}, \;\; i\in\{1,...,l\}{,\ i\ne s^\ast}
    \end{align*}
    {or }
    \begin{align*}
        {g(\x,\y^{(i)})-g(\x,\y^{(s^*)})=0,\;\;i=s^\ast.}
    \end{align*}
    Then, using the conditions in \eqref{eq:lagrange_mult}, the nonnegativity of the Lagrange multipliers $\mu_i^0,\mu_i^1 \geq 0, i=1,\ldots,k$, and the bound $\left|\frac{\partial \widetilde{f}(\x,\p)}{\partial p_i}\right| \leq L_f D, i\in\{1,...,k\}$ {from \eqref{eq:trivialresult}}, we obtain the following
        \begin{align*}
        g(\x,\y^{(i)})-g(\x,\y^{(s^*)})\leq\frac{\frac{\partial \widetilde{f}(\x,\p)}{\partial p_{s^*}}-\frac{\partial \widetilde{f}(\x,\p)}{\partial p_i}}{\gamma}\leq\frac{2L_f D}{\gamma}, \;\; i\in\{1,...,l\}{,\ i\ne s^\ast}
    \end{align*}
    {or}
    \begin{align*}
        {g(\x,\y^{(i)})-g(\x,\y^{(s^*)})=0\leq\frac{2L_f D}{\gamma}, \;\;i= s^\ast.}
    \end{align*}
    Next, following the same steps as in case 1, we conclude that 
    $$\sum\limits_{i=1}^kp_ig(\x,\y^{(i)})-\widetilde{V}_{\lambda}(\x)\leq\frac{5L_f D}{2\gamma}.$$
\end{proof}

\begin{proof}[\textbf{Proof of Proposition \ref{pro:disc_pen_optima}}]
\hfill \break
    Let $(\x^{\ast},\p^{\ast})$ be a KKT point of \eqref{eq:pen_reform}. By lemma \ref{lem:KKT1}, we obtain
    \begin{align}
        \sum\limits_{i=1}^{k} p^\ast_{i}g\left(\x^\ast,\y^{(i)}\right) - \widetilde{V}_{\lambda}(\x^\ast)\leq\frac{5L_f D}{2\gamma},
    \end{align}
    and set $\epsilon=\sum\limits_{i=1}^{k} p^{\ast}_{i}g\left(\x^{\ast},\y^{(i)}\right) - \widetilde{V}_{\lambda}(\x^{\ast})$. Also, there exist $(\mu_i^{0})^{\ast} \geq 0$, $(\mu_i^1)^{\ast} \geq 0$, $i=1,\ldots,k$, and $\mu^{\ast}$, such that the KKT conditions in \eqref{KKT:PEN} are satisfied.

    If we {define $G(\x,\p):=\sum\limits_{i=1}^{k} p_{i}g\left(\x,\y^{(i)}\right) - \widetilde{V}_{\lambda}(\x)$}, then the KKT conditions for the \eqref{eq:disc_reform} should satisfy the following system of equations:
    \begin{align}\label{eq:disc_kkt}
    \begin{pmatrix}
    \nabla_\x \widetilde{f}(\x,\p) \\
    \nabla_{p_1} \widetilde{f}(\x,\p)\\
    \vdots\\
    \nabla_{p_i} \widetilde{f}(\x,\p)\\
    \vdots\\
    \nabla_{p_k} \widetilde{f}(\x,\p)
    \end{pmatrix}
    +
    \lambda\begin{pmatrix}
    \nabla_\x G(\x,\p) \\
    g(\x,\y^{(1)})\\
    \vdots\\
    g(\x,\y^{(i)})\\
    \vdots\\
    g(\x,\y^{(k)})
    \end{pmatrix}
    -
    \sum\limits_{i=1}^{k}\begin{pmatrix}
    \mathbf{0} \\
    0\\
    \vdots\\
    \mu^0_i\\
    \vdots\\
    0
    \end{pmatrix}+\sum\limits_{i=1}^{k}\begin{pmatrix}
    \mathbf{0} \\
    0\\
    \vdots\\
    \mu^1_i\\
    \vdots\\
    0
    \end{pmatrix}
    +\begin{pmatrix}
    \mathbf{0} \\
    \mu\\
    \vdots\\
    \mu\\
    \vdots\\
    \mu
    \end{pmatrix}
    =0,
    \end{align}
    where $\mu_i^0 \geq 0$, $\mu_i^1 \geq 0$, $i=1,\ldots,k$, $\lambda\geq 0$ and $\mu$ are the Lagrange multipliers of the constraints ${\bf 0} \leq \p $, $\p \leq {\bf 1}$, $\sum\limits_{i=1}^{k} p_{i}g\left(\x,\y^{(i)}\right) - \widetilde{V}_{\lambda}(\x)\leq\epsilon$, and ${\bf 1}^\top \p=1$, respectively. 
    By comparing of the above conditions with the ones in \eqref{KKT:PEN}, it is not hard to see that a KKT point $(\x^{\ast},\p^{\ast})$ of \eqref{eq:pen_reform} satisfies the conditions in \eqref{eq:disc_kkt}, with $\mu=\mu^{\ast}$, $\mu_i^0=(\mu_i^0)^{\ast}$, $\mu_i^1=(\mu_i^1)^{\ast}$, $i=1,\ldots,k$ and $\lambda=\gamma$. It also satisfies the inequality $G(\x,\p) \leq \epsilon$ of \eqref{eq:disc_reform} for $\epsilon=G(\x^\ast,\p^\ast) < \frac{5L_f D}{2\gamma}$. Moreover, since we chose $\epsilon=G(\x^{\ast},\p^{\ast})$ it holds that $\lambda \left(G(\x^\ast,\p^\ast)-\epsilon \right)=0$, which means that the strict complementary condition holds for the $G(\x,\p) \leq \epsilon$ inequality. The satisfaction of the strict complementarity conditions regarding ${\bf 0} \leq \p $ and  $\p \leq {\bf 1}$ are guaranteed by the respective conditions already established for \eqref{eq:pen_reform}. Overall, a KKT point $(\x^{\ast},\p^{\ast})$ of \eqref{eq:pen_reform} is also a KKT point of \eqref{eq:disc_reform}.
    
    Furthermore, if $(\x^{\ast},\p^{\ast})$ is a local minimum of \eqref{eq:pen_reform}, then there exists a neighborhood $\mathcal{U}\subset \mathcal{X}\times\Delta$, such that for any $(\x,\p)$ in $U$, we have 
    \begin{align*}
        &f\left(\x,\sum\limits_{i=1}^{k} p_{i}\y^{(i)}\right) + \gamma \left[\sum\limits_{i=1}^{k} p_{i}g\left(\x,\y^{(i)}\right) - \widetilde{V}_{\lambda}(\x) \right] \\
        &\geq f\left(\x^{\ast},\sum\limits_{i=1}^{k} p^{\ast}_{i}\y^{(i)}\right) + \gamma \left[\sum\limits_{i=1}^{k} p^{\ast}_{i}g\left(\x^{\ast},\y^{(i)}\right) - \widetilde{V}_{\lambda}(\x^{\ast}) \right].
    \end{align*}
    Then for any $(\x,\p)$ in $\mathcal{U}$ such that $(\x,\p)$ is feasible for \eqref{eq:disc_reform}, we obtain
    \begin{align*}
        &f\left(\x^{\ast},\sum\limits_{i=1}^{k} p^{\ast}_{i}\y^{(i)}\right)\\
        &\leq f\left(\x,\sum\limits_{i=1}^{k} p_{i}\y^{(i)}\right)\\
        &\;\;+\gamma\left\{\left[\sum\limits_{i=1}^{k} p_{i}g\left(\x,\y^{(i)}\right) - \widetilde{V}_{\lambda}(\x) \right]-\left[\sum\limits_{i=1}^{k} p^{\ast}_{i}g\left(\x^{\ast},\y^{(i)}\right) - \widetilde{V}_{\lambda}(\x^{\ast}) \right]\right\}\\
        &\leq f\left(\x,\sum\limits_{i=1}^{k} p_{i}\y^{(i)}\right).
    \end{align*}
    {Since} any feasible point for \eqref{eq:disc_reform} satisfies the constraint
    \begin{align*}
        \sum\limits_{i=1}^{k} p_{i}g\left(\x,\y^{(i)}\right) - \widetilde{V}_{\lambda}(\x)\leq\epsilon,
    \end{align*}
    and here we selected $\epsilon=\sum\limits_{i=1}^{k} p^{\ast}_{i}g\left(\x^{\ast},\y^{(i)}\right) - \widetilde{V}_{\lambda}(\x^{\ast})${, then} $(\x^{\ast},\p^{\ast})$ is a local minimum of \eqref{eq:disc_reform}. The respective result for the global minima follows a reasoning similar to the above. 
\end{proof}

\subsection{Proofs of Section \ref{sec:alg}}\label{app:proofs:alg}

\begin{proof}[\textbf{Proof of Lemma \ref{lem:Lipschitz}}]
\hfill \break
Let us set $\widetilde{f}(\x,\p)=f\left(\x,\sum\limits_{i=1}^k p_i\y^{(i)}\right)$ and consider two arbitrary points $(\x_1,\p_1)$ {and} $(\x_2,\p_2)$ {belonging to} the set $\mathcal{X} \times \Delta$. Then, the following holds:
    \begin{align*}
        &\left\|\nabla_{\x,\p}\widetilde{f}(\x_1,\p_1)-\nabla_{\x,\p}\widetilde{f}(\x_2,\p_2)\right\|\\
        &=\left\|\nabla_{\x,\y}f(\x_1,\y_1)\nabla_{\x,\p}(\x_1,\y_1)-\nabla_{\x,\y}f(\x_2,\y_2)\nabla_{\x,\p}(\x_2,\y_2)\right\|\\
        &\leq\left\|\nabla_{\x,\y}f(\x_1,\y_1)\nabla_{\x,\p}(\x_1,\y_1)-\nabla_{\x,\y}f(\x_1,\y_1)\nabla_{\x,\p}(\x_2,\y_2)\right\|\\
        &\qquad +\left\|\nabla_{\x,\y}f(\x_1,\y_1)\nabla_{\x,\p}(\x_2,\y_2)-\nabla_{\x,\y}f(\x_2,\y_2)\nabla_{\x,\p}(\x_2,\y_2)\right\|\\
        &\leq\left\|\nabla_{\x,\y}f(\x_1,\y_1)\right\|\left\|\nabla_{\x,\p}(\x_1,\y_1)-\nabla_{\x,\p}(\x_2,\y_2)\right\|\\
        &\qquad +\left\|\nabla_{\x,\y}f(\x_1,\y_1)-\nabla_{\x,\y}f(\x_2,\y_2)\right\|\left\|\nabla_{\x,\p}(\x_2,\y_2)\right\|\\
        &\leq L_f\left\|\begin{pmatrix}
            \I-\I & \mathbf{0} & \cdots & \mathbf{0}\\
            \mathbf{0} & \y^{(1)}-\y^{(1)} & \cdots & \y^{(k)}-\y^{(k)}
        \end{pmatrix}\right\|_{SN}\\
        &\qquad +\overline{L}_f\left\|(\x_1,\p_1)-(\x_2,\p_2)\right\|\left\|\begin{pmatrix}
            \I & \mathbf{0} & \cdots & \mathbf{0}\\
            \mathbf{0} & \y^{(1)} & \cdots & \y^{(k)}
        \end{pmatrix}\right\|_{SN}\\
        &=\overline{L}_f\left\|(\x_1,\p_1)-(\x_2,\p_2)\right\|\left\|\begin{pmatrix}
            \I & \mathbf{0} & \cdots & \mathbf{0}\\
            \mathbf{0} & \y^{(1)} & \cdots & \y^{(k)}
        \end{pmatrix}\right\|_{SN}
    \end{align*}
    where $\y_{1}=\sum\limits_{i=1}^k p_{1_{i}}\y^{(i)}$, $\y_{2}=\sum\limits_{i=1}^k p_{2_{i}}\y^{(i)}$, with $p_{1_{i}},p_{2_{i}}$ the $i$th elements of $\p_1,\p_2$, respectively, and $\|\cdot\|_{SN}$ denotes the {spectral} norm. Note that
    \begin{align*}
        \begin{pmatrix}
            \I & \mathbf{0} & \cdots & \mathbf{0}\\
            \mathbf{0} & \y^{(1)} & \cdots & \y^{(k)}
        \end{pmatrix}\begin{pmatrix}
            \I & \mathbf{0} & \cdots & \mathbf{0}\\
            \mathbf{0} & \y^{(1)} & \cdots & \y^{(k)}
        \end{pmatrix}^\top
        =\begin{pmatrix}
            \I & \mathbf{0} \\
            \mathbf{0} & \left\|\y^{(1)}\right\|^2+\cdots+\left\|\y^{(k)}\right\|^2
        \end{pmatrix},
    \end{align*}
    so
    \begin{align*}
        \left\|\begin{pmatrix}
            \I & \mathbf{0} & \cdots & \mathbf{0}\\
            \mathbf{0} & \y^{(1)} & \cdots & \y^{(k)}
        \end{pmatrix}\right\|_{SN}\leq\max\{1,\sqrt{k}D^2\}.
    \end{align*}
    Consequently, 
    \begin{align}\label{eq:grad_xp_ftilde}
        \left\|\nabla_{\x,\p}\widetilde{f}(\x_1,\p_1)-\nabla_{\x,\p}\widetilde{f}(\x_2,\p_2)\right\|\leq\overline{L}_f\max\{1,\sqrt{k}D^2\}\left\|(\x_1,\p_1)-(\x_2,\p_2)\right\|.
    \end{align}

    Furthermore, we set $\widetilde{g}(\x,\p)=\sum_{i=1}^k p_ig(\x,\y^{(i)})$ and consider two arbitrary points $(\x_1,\p_1)$, $(\x_2,\p_2)$ of the set $\mathcal{X} \times \Delta$. 
    \begin{align}\label{eq:grad_xp_gtilde}
         &\left\|\nabla_{\x,\p}\widetilde{g}(\x_1,\p_1)-\nabla_{\x,\p}\widetilde{g}(\x_2,\p_2)\right\| \nonumber\\
        =&\left\|\begin{pmatrix}
            \sum_{i=1}^k\left({p_1}_i\nabla_{\x}g(\x_1,\y^{(i)})-{p_2}_i\nabla_{\x}g(\x_2,\y^{(i)})\right)\\
            g(\x_1,\y^{(1)})-g(\x_2,\y^{(1)})\\
            \vdots \\
            g(\x_1,\y^{(k)})-g(\x_2,\y^{(k)})\\
        \end{pmatrix}\right\|\nonumber\\
        =&\left\|\begin{pmatrix}
            \sum_{i=1}^k\left({p_1}_i(\nabla_{\x}g(\x_1,\y^{(i)})-\nabla_{\x}g(\x_2,\y^{(i)}))+({p_1}_{i}-{p_2}_i)\nabla_{\x}g(\x_2,\y^{(i)})\right)\\
            g(\x_1,\y^{(1)})-g(\x_2,\y^{(1)})\\
            \vdots\\
            g(\x_1,\y^{(k)})-g(\x_2,\y^{(k)})\\
        \end{pmatrix}\right\|\nonumber\\
        \leq&\Biggl[\left(\sum_{i=1}^k\left({p_1}_i(\nabla_{\x}g(\x_1,\y^{(i)})-\nabla_{\x}g(\x_2,\y^{(i)}))+({p_1}_{i}-{p_2}_i)\nabla_{\x}g(\x_2,\y^{(i)})\right)\right)^2\nonumber\\
        &+kL_g^2\left\|\x_1-\x_2\right\|^2\Biggr]^{\frac{1}{2}}\nonumber\\
        \leq&\left(\left(\overline{L}_g\left\|\x_1-\x_2\right\|+L_g\left\|\p_1-\p_2\right\|\right)^2+kL_g^2\left\|\x_1-\x_2\right\|^2\right)^{\frac{1}{2}}\nonumber\\
        \leq&\left(\left(\overline{L}_g+L_g\right)^2+kL_g^2\right)^{\frac{1}{2}}\left\|(\x_1,\p_1)-(\x_2,\p_2)\right\|
    \end{align}
    where in the first inequality we used Assumption \ref{ass:basics}\ref{ass:basics:Lip_bound_g} and in the second inequality we used Assumptions \ref{ass:basics}\ref{ass:basics:Lip_grad_g} and \ref{ass:basics}\ref{ass:basics:Lip_bound_g}.
    
    Recall that, by {Proposition} \ref{pro:dvf_properties}, $\nabla \widetilde{V}_{\lambda}(x)$ is $\left(\frac{L_g^2k}{\lambda}+\overline{L}_g\right)$-Lipschitz continuous. By combining \eqref{eq:grad_xp_ftilde}, \eqref{eq:grad_xp_gtilde} and the Lipschitz continuity of $\nabla \widetilde{V}_{\lambda}(x)$ we get the following  
    \begin{align*}
        &\left\|\nabla_{\x,\p}F_\gamma(\x_1,\p_1)-\nabla_{\x,\p}F_\gamma(\x_2,\p_2)\right\|\\
        =&\bigg\|\left(\nabla_{\x,\p}\widetilde{f}(\x_1,\p_1)-\nabla_{\x,\p}\widetilde{f}(\x_2,\p_2)\right)+\gamma\bigg[\left(\nabla_{\x,\p}\widetilde{g}(\x_1,\p_1)-\nabla_{\x,\p}\widetilde{g}(\x_2,\p_2)\right)\\
        &\hspace{5.4cm}+\left(\nabla_{\x,\p}\widetilde{V}_{\lambda}(x_1)-\nabla_{\x,\p}\widetilde{V}_{\lambda}(x_2)\right)\bigg]\bigg\|\nonumber\\
        \leq&\Biggl\{\overline{L}_f\max\{1,\sqrt{k}D^2\}+ \\ &+\gamma\left[\left(\left(\overline{L}_g+L_g\right)^2+kL_g^2\right)^{\frac{1}{2}}+\frac{L_g^2k}{\lambda}+\overline{L}_g\right]\Biggr\}\left\|(\x_1,\p_1)-(\x_2,\p_2)\right\|.
    \end{align*}
\end{proof}

\begin{proof}[\textbf{Proof of Theorem \ref{thm:convergence}}]
\hfill \break
    This theorem consists of two parts. Below we provide the proof of each part.
    
    1. Let $\z=(\x,\p)$. By Lipschitz smoothness of $F_{\gamma}(\z)$ and under the assumption that $\alpha \leq L_{\gamma}^{-1}$, it holds that
    \begin{align}\label{eq:conv_rate_descent_lemma}
        F_{\gamma}(\z_{t+1})&\leq F_{\gamma}(\z_t)+\left<\nabla F_{\gamma}(\z_t),\z_{t+1}-\z_t\right>+\frac{L_\gamma}{2}\left\|\z_{t+1}-\z_t\right\|^2 \nonumber\\
        &\leq F_{\gamma}(\z_t)+\left<\nabla F_{\gamma}(\z_t),\z_{t+1}-\z_t\right>+\frac{1}{2\alpha}\left\|\z_{t+1}-\z_t\right\|^2.%\label{eq:conv_rate_descent_lemma}
    \end{align}
    Note that $\z_{t+1}=\text{Proj}_{\mathcal{X}\times\Delta}(\z_t-\alpha \nabla F_{\gamma}(\z_t))$. Then, with the use of the optimality condition of the projection operator we obtain the following
    \begin{align*}
        \left<\z_{t+1}-\z,\z_{t+1}-(\z_t-\alpha\nabla F_{\gamma}(\z_t))\right>\leq 0,
    \end{align*}
    which holds for every $\z \in \mathcal{X}\times\Delta$. Particularly we can select $\z$ as $\z_t$ and obtain
    \begin{align*}
        \left<\z_{t+1}-\z_t,\z_{t+1}-(\z_t-\alpha\nabla F_{\gamma}(\z_t))\right>\leq 0,
    \end{align*}
    which means
    \begin{align}\label{eq:inner_product}
        \left<\nabla F_\gamma(\z_k),\z_{t+1}-\z_t\right>\leq-\frac{1}{\alpha}\left\|\z_{t+1}-\z_t\right\|^2.
    \end{align}
    Substituting (\ref{eq:inner_product}) into (\ref{eq:conv_rate_descent_lemma}), it  follows that
    \begin{align*}
        \left\|\frac{1}{\alpha}\left(\z_{t+1}-\z_t\right)\right\|^2\leq 2\frac{F_{\gamma}(\z_t)-F_{\gamma}(\z_{t+1})}{\alpha}
    \end{align*}
    Given that $\widetilde{f}(\z) \geq C_f$ and $\widetilde{g}(\z) - \widetilde{V}_{\lambda}(\x) \geq -\frac{\lambda}{2}$ for all $\z \in \mathcal{X} \times \Delta$, it follows that $F{\gamma}(\z)$ is bounded below by $C_f-\frac{\lambda}{2}$ across this domain. To demonstrate convergence, consider:
    \begin{align}
        \sum^T_{t=1}\left\|\frac{1}{\alpha}\left(\z_{t+1}-\z_t\right)\right\|^2\leq 2\frac{F_\gamma(\z_1)-C_f+\frac{\lambda}{2}}{\alpha}.
    \end{align}
    By multiplying both sides with $\frac{1}{T}$, we prove the convergence.
    
    2. Suppose $\lim_{t\to\infty}\z_t=\z^{\ast}$. It is clear that $\z^{\ast}$ is a fixed point of the algorithm, which means $\z^{\ast}=\text{Proj}_{\mathcal{X}\times\Delta}\left(\z^{\ast}-\alpha\nabla F_{\gamma}(\z^{\ast})\right)$. According to the optimality condition of the projection operator, the following holds 
    \begin{align*}
        \left<\z^\ast-\left(\z^\ast-\alpha\nabla F_\gamma(\z^{\ast})\right),\z^{\ast}-\z\right>\leq 0,\ \forall \z\in\mathcal{X}\times\Delta
    \end{align*}
    Since $\alpha>0$, simplifying this expression, we obtain
    \begin{align*}
        \left<\nabla F_\gamma(\z^{\ast}),\z^{\ast}-\z\right>\leq 0,\ \forall \z\in\mathcal{X}\times\Delta
    \end{align*}
    which indicates that $\z^{\ast}$ is a stationary point.
\end{proof}

\clearpage
\addcontentsline{toc}{section}{References}

\bibliographystyle{unsrtnat}
\bibliography{refs}

\end{document}